\documentclass[12pt,a4paper,fleqn]{article}
\usepackage{a4wide}
\usepackage{amsfonts}
\usepackage{amsmath}
\usepackage{latexsym}
\usepackage{amssymb}
\usepackage{euscript}
\usepackage{graphicx}
\usepackage{units}
\usepackage{mathrsfs}

\usepackage{graphicx}
\usepackage{color}
\usepackage{amssymb}
\usepackage{amssymb}
\usepackage[T1]{fontenc}
\usepackage{latexsym}
\usepackage{xypic}
\usepackage{eufrak}
\usepackage{euscript}
\usepackage{amsfonts,amsmath}
\usepackage{verbatim}
\usepackage{fancyhdr}
\usepackage[english]{babel}

\newtheorem{prop}{Proposition}[section]

\newtheorem{corollary}[prop]{Corollary}
\newtheorem{proposition}[prop]{Proposition}
\newtheorem{lemme}[prop]{Lemma}
\newtheorem{rem}[prop]{Remark}
\newtheorem{remark}[prop]{Remark}

\newtheorem{thm}[prop]{Theorem}
\newtheorem{theorem}[prop]{Theorem}

\renewcommand{\geq}{\geqslant}
\def\leq{\leqslant}

\newcommand{\N}{\mathbb{N}}
\newcommand{\Z}{\mathbb{Z}}
\newcommand{\R}{\mathbb{R}}

\def\cal{\mathcal}
\def\1{{\mathbf{1}}}

\def\1{{\mathbf{1}}}
\def\0.5{{\frac{1}{2}}}

\newcommand{\fin}
{ \vspace{-0.6cm}
\begin{flushright}
\mbox{$\Box$}
\end{flushright}
\noindent }

\newenvironment{proof}[1]{\begin{trivlist}\item {\it
\bf Proof.}\quad} {\qed\end{trivlist}}

\newcommand{\qed}{\nopagebreak\hspace*{\fill}
{\vrule width6pt height6ptdepth0pt}\par}

\begin{document}

\begin{center}
{\large{\bf Weighted power variations
of iterated Brownian motion}}\\~\\
by Ivan Nourdin\footnote{Laboratoire de Probabilit{\'e}s et
Mod{\`e}les Al{\'e}atoires, Universit{\'e} Pierre et Marie Curie,
Bo{\^\i}te courrier 188, 4 Place Jussieu, 75252 Paris Cedex 5,
France, {\tt ivan.nourdin@upmc.fr}} and Giovanni
Peccati\footnote{ Laboratoire de Statistique Th\'eorique et
Appliqu\'ee, Universit\'e Pierre et Marie Curie, 8\`eme \'etage,
b\^atiment A, 175 rue du Chevaleret, 75013 Paris, France,
{\tt giovanni.peccati@gmail.com}}\\
{\it Universit\'e Paris VI}\\~\\
\end{center}
{\small \noindent {\bf Abstract:} We characterize the asymptotic
behaviour of the weighted power variation processes associated
with iterated Brownian motion. We prove weak convergence results
in the sense of finite dimensional distributions, and show that
the laws of the limiting objects can always be expressed in terms
of three independent Brownian motions $X$, $Y$ and $B$, as well as
of the local times of $Y$. In particular, our results involve
``weighted'' versions of Kesten and Spitzer's \textsl{Brownian
motion in random scenery}. Our findings extend the
theory initiated by Khoshnevisan and Lewis (1999), and should be
compared with the recent result by Nourdin and R\'eveillac
(2008), concerning the weighted power
variations of fractional Brownian motion with Hurst index $H=1/4$.\\


\\

\section{Introduction and main results}

The characterization of the single-path behaviour of a given
stochastic process is often based on the study of its power
variations. A quite extensive literature has been developed on the
subject, see e.g. \cite{CorNuWoe, NNT} (as well as the forthcoming
discussion) for references concerning the power variations of
Gaussian and Gaussian-related processes, and \cite{BaNiGrSh} (and
the references therein) for applications of power variation
techniques to the continuous-time modeling of financial markets.
Recall that, for a given real $\kappa>1$ and a given real-valued
stochastic process $Z$, the $\kappa$-power variation of $Z$, with
respect to a partition $\pi=\{0= t_{0}<t_{1}<\ldots<t_{N}=1\}$ of
$[0,1]$ ($N\geq 2$ is some integer),
is defined to be the sum
\begin{equation}\label{SUM_absolute}
\sum_{k=1}^{N}|Z_{t_{k}}-Z_{t_{k-1}}|^{\kappa}.
\end{equation}
For the sake of simplicity, from now on we shall only consider the
case where $\pi$ is a dyadic partition, that is, $N = 2^{n}$ and
$t_{k}=k2^{-n}$, for some integer $n\geq 2$ and for
$k\in\{0,\ldots,2^n\}$.

The aim of this paper is to study the asymptotic behaviour, for
every integer $\kappa\geq 2$ and for $n\rightarrow\infty$, of the
(dyadic) $\kappa$-power variations associated with a remarkable
non-Gaussian and self-similar process with stationary increments,
known as \textsl{iterated Brownian motion} (in the sequel,
I.B.M.). Formal definitions are given below: here, we shall only
observe that I.B.M. is a self-similar process of order $\frac14$,
realized as the composition of two independent Brownian motions.
As such, I.B.M. can be seen as a non-Gaussian counterpart to
Gaussian processes with the same order of self-similarity, whose
power variations (and related functionals) have been recently the
object of an intense study. In this respect, the study of the
single-path behaviour of I.B.M. is specifically relevant, when one
considers functionals that are obtained from (\ref{SUM_absolute})
by dropping the absolute value (when $\kappa$ is odd), and by
introducing some weights. More precisely, in what follows we shall
focus on the asymptotic behaviour of \textsl{weighted variations}
of the type
\begin{eqnarray}\label{sum-intro}
&&\sum_{k=1}^{2^n}f(Z_{(k-1)2^{-n}})\big(Z_{k2^{-n}}-Z_{(k-1)2^{-n}}\big)^\kappa,\quad\kappa=2,3,4,\ldots,
\end{eqnarray}
or
\begin{eqnarray}\label{sum-intro2}
&&\sum_{k=1}^{2^n}\frac12[f(Z_{(k-1)2^{-n}})+f(Z_{k2^{-n}})]
\big(Z_{k2^{-n}}-Z_{(k-1)2^{-n}}\big)^\kappa,\quad\kappa=2,3,4,\ldots,
\end{eqnarray}
for a real function $f:\R\rightarrow\R$ satisfying some suitable
regularity conditions.

Before dwelling on I.B.M., let us recall some recent results
concerning (\ref{sum-intro}), when $Z=B$ is a \textsl{fractional Brownian
motion} (fBm) of Hurst index $H \in (0,1)$ (see e.g. \cite{Nbook} for
definitions) and, for instance, $\kappa\geq 2$ is an {\sl even} integer.
Recall that, in particular, $B$ is a continuous Gaussian process,
with an order of self-similarity equal to $H$. In what follows,
$f$ denotes a smooth enough function such that $f$ and its
derivatives have subexponential growth. Also, here and for the
rest of the paper, $\mu_q$, $q\geq1$, stands for the $q$th moment
of a standard Gaussian random variable, that is, $\mu_q =0$ if $q$
is odd, and
\begin{equation}\label{mu} \mu_q =
\frac{q!}{2^{q/2}(q/2)!},\quad\mbox{if $q$ is even.}
\end{equation}
We have (see \cite{Nourdin,NNT,NR}) as $n\rightarrow\infty$:
\begin{enumerate}
\item When $H>\frac34$,
\begin{equation}\label{NNT1}
2^{n-2Hn} \sum_{k=1}^{2^n}\!f(B_{(k-1)2^{-n}})\! \left[\big(2^{n
H}(B_{k2^{-n}}\!-\!B_{(k-1)2^{-n}})\big)^\kappa\!\!-\!\!\mu_{\kappa}\right]
\!\!\,\,\,{\stackrel{{\rm L}^2}{\longrightarrow}}\,\,\,\!\!\!
\mu_{\kappa-2}\binom{\kappa}{2}\int_0^1 f(B_s)dZ^{(2)}_s,
\end{equation}
where $Z^{(2)}$ denotes the Rosenblatt process canonically constructed 
from $B$ (see \cite{NNT} for more details.)
\item When $H = \frac34$,
\begin{equation}\label{NNT2}
\frac{2^{-\frac{n}{2}}}{\sqrt{n}}\sum_{k=1}^{2^n}\!f(B_{(k-1)2^{-n}})\!
\left[\big(2^{\frac{3n}{4}}(B_{k2^{-n}}\!\!-\!\!B_{(k-1)2^{-n}})\big)^\kappa\!-\!\mu_{\kappa}\right]
\!\!\,\,\,{\stackrel{{\rm Law}}{\longrightarrow}}\,\,\,\!\!
\sigma_{\frac34,\kappa}\int_0^1 f(B_s)dW_s,
\end{equation}
where $\sigma_{\frac34,\kappa}$ is an explicit constant depending only on
$\kappa$, and $W$ is a standard Brownian motion independent of
$B$.
\item When $\frac14<H<\frac34$,
\begin{equation}\label{NNT3}
2^{-\frac{n}{2}}\sum_{k=1}^{2^n}\!f(B_{(k-1)2^{-n}})\!
\left[\big(2^{n
H}(B_{k2^{-n}}\!-\!B_{(k-1)2^{-n}})\big)^\kappa\!-\!\mu_{\kappa}\right]
\!\!\,\,\,{\stackrel{{\rm Law}}{\longrightarrow}}\,\,\,\!\!
\sigma_{H,\kappa}\!\!\int_0^1 \!\!f(B_s)dW_s,
\end{equation}
for an explicit constant $\sigma_{H,\kappa}$ depending only on $H$
and $\kappa$, and where $W$ denotes a standard Brownian motion
independent of $B$.
\item When $H=\frac14$,
\begin{eqnarray}
2^{-\frac{n}2}\! \sum_{k=1}^{2^n}\!f(B_{(k-1)2^{-n}}) \!\left[\big(2^{n
H}(B_{k2^{-n}}\!-\!B_{(k-1)2^{-n}})\big)^\kappa\!-\!\mu_{\kappa}\right]\!\!
\,\,\,{\stackrel{{\rm Law}}{\longrightarrow}}\,\,\,\!\!\frac14\mu_{\kappa-2}\binom{\kappa}{2}
\! \int_0^1\! f''(B_s)ds \notag\\
+ \sigma_{\frac14,\kappa} \int_0^1 f(B_s)dW_s,\quad\label{NouRev}
\end{eqnarray}
where $\sigma_{\frac14,\kappa}$ is an explicit constant depending only on
$\kappa$, and $W$ is a standard Brownian motion independent of
$B$.
\item When $H<\frac14$,
\begin{equation}\label{NNT4}
2^{2Hn-n}\! \sum_{k=1}^{2^n}\!f(B_{(k-1)2^{-n}}) \!\left[\big(2^{n
H}(B_{k2^{-n}}\!-\!B_{(k-1)2^{-n}})\big)^\kappa\!-\!\mu_{\kappa}\right]\!\!
\,\,\,{\stackrel{{\rm
L}^2}{\longrightarrow}}\,\,\,\!\!\frac14\mu_{\kappa-2}\binom{\kappa}{2}
\! \int_0^1\! f''(B_s)ds.
\end{equation}
\end{enumerate}

In the current paper, we focus on the iterated Brownian motion, 
which is a continuous non-Gaussian self-similar process of order
$\frac14$. More precisely, let $X$ be a two-sided Brownian motion, and let $Y$ be a standard
(one-sided) Brownian motion independent of $X$. In what follows,
we shall denote by $Z$ the {\sl iterated Brownian motion}
\textit{(I.B.M.) associated with} $X$ \textit{and} $Y$, that is,
\begin{equation}
Z(t)=X\big(Y(t)\big),\quad t\geq 0. \label{IBM!}
\end{equation}
The process $Z$ appearing in (\ref{IBM!}) has been first (to the best of our knowledge)
introduced in \cite{Burdzy}, 
and then further studied in a number
of papers -- see for instance \cite{KL2} for a comprehensive
account up to 1999, and \cite{deBlassie, Nane1, Nane2, OB} for
more recent references on the subject. Such a process can be
regarded as the realization of a Brownian motion on a random
fractal (represented by the path of the underlying motion $Y$).
Note that $Z$ is self-similar of order $\frac14$, $Z$ has
stationary increments, and $Z$ is neither a Dirichlet process nor
a semimartingale or a Markov process in its own filtration. A
crucial question is therefore how one can define a stochastic
calculus with respect to $Z$. This issue has been tackled by
Khoshnevisan and Lewis in the ground-breaking paper \cite{KL1}
(see also \cite{KL2}), where the authors develop a
Stratonovich-type stochastic calculus with respect to $Z$, by
extensively using techniques based on the properties of some
special arrays of Brownian stopping times, as well as on
excursion-theoretic arguments. Khoshnevisan and Lewis' approach
can be roughly summarized as follows. Since the paths of $Z$ are
too irregular, one cannot hope to effectively define stochastic
integrals as limits of Riemann sums with respect to a
deterministic partition of the time axis. However, a winning idea
is to approach deterministic partitions by means of random
partitions defined in terms of hitting times of the underlying
Brownian motion $Y$. In this way, one can bypass the random
``time-deformation'' forced by (\ref{IBM!}), and perform
asymptotic procedures by separating the roles of $X$ and $Y$ in
the overall definition of $Z$. Later in this section, by adopting
the same terminology introduced in \cite{KL2}, we will show that
the role of $Y$ is specifically encoded by the so-called
``intrinsic skeletal structure'' of $Z$.

By inspection of the techniques developed in \cite{KL1}, one sees
that a central role in the definition of a stochastic calculus
with respect to $Z$ is played by the asymptotic behavior of the
quadratic, cubic and quartic variations associated with $Z$. Our
aim in this paper is to complete the results of \cite{KL2}, by
proving asymptotic results involving \textit{weighted} power
variations of $Z$ of arbitrary order, where the weighting is
realized by means of a well-chosen real-valued function of $Z$.
Our techniques involve some new results concerning the weak
convergence of non-linear functionals of Gaussian processes,
recently proved in \cite{PT}. As explained above, our results
should be compared with the recent findings, concerning power
variations of Gaussian processes, contained in \cite{Nourdin, NNT,
NR}. 

Following Khoshnevisan and Lewis \cite{KL1,KL2}, we start by
introducing the so-called {\sl intrinsic skeletal structure} of
the I.B.M. $Z$ appearing in (\ref{IBM!}). This structure is
defined through a sequence of collections of stopping times (with
respect to the natural filtration of $Y$), noted
\begin{equation}
\mathscr{T}_n=\{T_{k,n}: k\geq0\}, \quad n\geq1, \label{TN}
\end{equation}
which are in turn expressed in terms of the subsequent hitting
times of a dyadic grid cast on the real axis. More precisely, let
$\mathscr{D}_n= \{j2^{-n/2}:\,j\in\Z\}$, $n\geq 1$, be the dyadic
partition (of $\R$) of order $n/2$. For every $n\geq 1$, the
stopping times $T_{k,n}$, appearing in (\ref{TN}), are given by
the following recursive definition: $T_{0,n}= 0$, and
$$
T_{k,n}= \inf\big\{s>T_{k-1,n}:\quad
Y(s)\in\mathscr{D}_n\setminus\{Y(T_{k-1,n})\}\big\},\quad k\geq 1,
$$
where, as usual, $A\setminus B = A\cap B^c$ ($B^c$ is the
complement of $B$). Note that the definition of $T_{k,n}$, and
therefore of $\mathscr{T}_n$, only involves the one-sided Brownian
motion $Y$, and that, for every $n\geq1$, the discrete stochastic
process
$$\mathscr{Y}_n=\{Y(T_{k,n}):k\geq0\}$$ defines a simple random
walk over $\mathscr{D}_n$. The intrinsic skeletal structure of $Z$
is then defined to be the sequence $${\rm I.S.S.} =
\{\mathscr{D}_n,\mathscr{T}_n,\mathscr{Y}_n: n\geq1 \},$$
describing the random scattering of the paths of $Y$ about the
points of the partitions $\{\mathscr{D}_n\}$. As shown in
\cite{KL1}, the I.S.S. of $Z$ provides an appropriate sequence of
(random) partitions upon which one can build a stochastic calculus
with respect to $Z$. It can be shown that, as $n$ tends to
infinity, the collection $\{T_{k,n}:\,k\geq 0\}$ approximates the
common dyadic partition $\{k2^{-n}:\,k\geq 0\}$ of order $n$ (see
\cite[Lemma 2.2]{KL1} for a precise statement). Inspired again by
\cite{KL1}, we shall use the I.S.S. of $Z$ in order to define and
study weighted power variations, which are the main object of this
paper. To this end, recall that $\mu_\kappa$ is defined, via
(\ref{mu}), as the $\kappa$th moment of a centered standard
Gaussian random variable. Then, the {\sl weighted power variation}
of the I.B.M. $Z$, associated with a real-valued function $f$,
with an instant $t\in[0,1]$, and with integers $n\geq 1$ and
$\kappa\geq2$, is defined as follows:
\begin{eqnarray}
V_n^{(\kappa)}(f,t)&=&\!\!\!\frac12\!\!\sum_{k=1}^{\lfloor 2^n t\rfloor}
\!\!\left(
f\big(Z(T_{k,n})\big) + f\big(Z(T_{k-1,n})\big)\right)\left(
\big(
Z(T_{k,n})-Z(T_{k-1,n})
\big)^\kappa - \mu_\kappa\,2^{-\kappa\frac{n}4}
\right).\quad\mbox{ }\quad\mbox{ }\label{pow2}
\end{eqnarray}
Note that, due to self-similarity and independence,
$$\mu_\kappa\,2^{-\kappa\frac{n}4} = E\big[\big(
Z(T_{k,n})-Z(T_{k-1,n}) \big)^\kappa\big]=E\big[\big(
Z(T_{k,n})-Z(T_{k-1,n}) \big)^\kappa \mid Y \big ].$$ For each
integer $n\geq 1$, $k\in\Z$ and $t\geq 0$, let $U_{j,n}(t)$ (resp.
$D_{j,n}(t)$) denote the number of \textit{upcrossings} (resp.
\textit{downcrossings}) of the interval
$[j2^{-n/2},(j+1)2^{-n/2}]$ within the first $\lfloor 2^n
t\rfloor$ steps of the random walk $\{Y(T_{k,n})\}_{k\geq 1}$ (see
formulae (\ref{UPC}) and (\ref{DOWNC}) below for precise
definitions). The following lemma plays a crucial role in the
study of the asymptotic behavior of $V_n^{(\kappa)}(f,\cdot)$:
\begin{lemme}\label{lm-kl} (See \cite[Lemma 2.4]{KL1}) Fix $t\in[0,1]$, $\kappa\geq 2$
and let $f:\R\rightarrow\R$ be any real function.
Then
\begin{eqnarray}
&&V_n^{(\kappa)}(f,t)= \frac12\sum_{j\in\Z}
\left(
f(X_{\left( j-1\right) 2^{-\frac{n}2}})
+f(X_{j\, 2^{-\frac{n}2}})
\right)\times \label{amenable}\\
&&\hskip3cm
 \left[ \left( X_{j2^{-\frac{n}2}}-
X_{(j-1)2^{-\frac{n}2}}\right)^\kappa
-\mu_\kappa\,2^{-\kappa\frac{n}4} \right]
\big(U_{j,n}(t)+(-1)^\kappa\,D_{j,n}(t)\big). \notag
\end{eqnarray}
\end{lemme}
The main feature of the decomposition (\ref{amenable}) is that it
separates $X$ from $Y$, providing a representation of
$V_n^{(\kappa)}(f,t)$ which is amenable to analysis. Using Lemma
\ref{lm-kl} as a key ingredient, Khoshnevisan and Lewis \cite{KL1}
proved the following results, corresponding to the case where $f$
is identically one in (\ref{pow2}): as $n\rightarrow\infty$,
$$
\frac{2^{-n/4}}{\sqrt{2}}\,\, V_n^{(2)}(1,\cdot)
\,\,\,\overset{D[0,1]}{\Longrightarrow}\,\,\,
{\rm B.M.R.S.}
\quad\mbox{and}\quad
\frac{2^{n/4}}{\sqrt{96}}\,\, V_n^{(4)}(1,\cdot)
\,\,\,\overset{D[0,1]}{\Longrightarrow}\,\,\,
{\rm B.M.R.S.}
$$
Here, and for the rest of the paper,
$\overset{D[0,1]}{\Longrightarrow}$ stands for the convergence in
distribution in the Skorohod space $D[0,1]$, while `B.M.R.S.'
indicates Kesten and Spitzer's {\sl Brownian Motion in Random
Scenery} (see \cite{KS}). This object is defined as:
\begin{equation}
{\rm B.M.R.S.}=\left\{\int_\R L^x_t(Y)dB_x\right\}_{t\in[0,1]},
\label{BMRS}
\end{equation}
where $B$ is a two-sided Brownian motion independent of $X$ and
$Y$, and $\{L^x_t(Y)\}_{x\in\R,\,t\in[0,1]}$ is a jointly
continuous version of the local time process of $Y$ (the
independence of $X$ and $B$ is immaterial here, and will be used
in the subsequent discussion). In \cite{KL1} it is also proved
that the asymptotic behavior of the cubic variation of $Z$ is very
different, and that in this case the limit is I.B.M. itself,
namely:
$$
\frac{2^{n/2}}{\sqrt{15}}\,\,V_n^{(3)}(1,\cdot)
\,\,\,\overset{D[0,1]}{\Longrightarrow}\,\,\,
{\rm I.B.M.}
$$

As anticipated, our aim in the present paper is to characterize
the asymptotic behavior of $V_n^{(\kappa)}(f,t)$ in (\ref{pow2}),
$n\rightarrow\infty$, in the case of a general function $f$ and of
a general integer $\kappa\geq2$. Our main result is the following:
\begin{thm}\label{thm1}
Let $f:\R\rightarrow\R$ belong to ${\rm C}^2$ with $f'$ and $f''$ bounded, and
$\kappa\geq 2$ be an integer.
Then, as $n\rightarrow\infty$,
\begin{enumerate}
\item if $\kappa$ is even,
$\left\{X_x,2^{(\kappa-3)\frac{n}4}\,V_n^{(\kappa)}(f,t)\right\}_{x\in\R,\,t\in[0,1]}$
converges in the sense of finite dimensional distributions
(f.d.d.) to
\begin{eqnarray}
\left\{X_x,\sqrt{\mu_{2\kappa}-\mu_\kappa^2}\int_\R f(X_z)L_t^z(Y)dB_z\right\}_{x\in\R,\,t\in[0,1]};\label{gen2}
\end{eqnarray}
\item if $\kappa$ is odd,
$\left\{X_x,
2^{(\kappa-1)\frac{n}4}\,V_n^{(\kappa)}(f,t)\right\}_{x\in\R,\,t\in[0,1]}$
converges in the sense of f.d.d. to
\begin{eqnarray}
\left\{X_x,\int_0^{Y_t} f(X_z)\big(\mu_{\kappa+1}\,d^\circ X_z
+\sqrt{\mu_{2\kappa}-\mu_{\kappa+1}^2}\,dB_z\big)\right\}_{x\in\R,\,t\in[0,1]}.\label{gen3}
\end{eqnarray}
\end{enumerate}
\end{thm}
\begin{rem}\label{rk13}
{\rm
\begin{enumerate}
\item Here, and for the rest of the paper, we shall write
$\int_0^\cdot f(X_z)d^\circ X_z$ to indicate the {\sl Stratonovich} integral of $f(X)$ with respect to
the Brownian motion $X$. On the other hand, $\int_0^t f(X_z)dB_z$ (resp. $\int_0^t f(X_z)L_t^z(Y)dB_z$) 
is well-defined (for each $t$) as the {\sl Wiener-It\^{o}} stochastic integral of the random mapping
$z\mapsto f(X_z)$ (resp. $z\mapsto f(X_z)L_t^z(Y)$), with respect to the {\sl independent}
Brownian motion $B$. In particular, one uses the fact that the
mapping $z\mapsto L_t^z(Y)$ has a.s. compact support.
\item We call the process
\begin{equation}\label{WBMRS}
t\mapsto \int_\R f(X_z)L_t^z(Y)dB_z
\end{equation}
appearing in (\ref{gen2}) a \textsl{Weighted Browian Motion in
Random Scenery} (W.B.M.R.S. -- compare with (\ref{BMRS})), the
weighting being randomly determined by $f$ and by the independent
Brownian motion $X$. 
\item The relations (\ref{gen2})-(\ref{gen3}) can be reformulated in the sense of ``stable convergence''.
For instance, (\ref{gen2}) can be rephrased by saying that the
finite dimensional distributions of
$$2^{(\kappa-3)\frac{n}4}\,V_n^{(\kappa)}(f,\cdot)$$
converge {\sl $\sigma(X)$-stably} to those of
 $$\sqrt{\mu_{2\kappa}-\mu_\kappa^2}\int_\R
f(X_z)L^z_\cdot(Y)dB_z$$(see e.g. Jacod and Shiryayev \cite{JS}
for an exhaustive discussion of stable convergence).
\item Of course, one recovers finite dimensional versions of the results by Khoshnevisan and Lewis by choosing
$f$ to be identically one in (\ref{gen2})-(\ref{gen3}).
\item To keep the length of this paper within bounds, we defer to future analysis the
rather technical investigation of the tightness of the processes $
2^{(\kappa-3)\frac{n}4}\,V_n^{(\kappa)}(f,t)$ ($\kappa$ even) and
$2^{(\kappa-1)\frac{n}4}\,V_n^{(\kappa)}(f,t)$ ($\kappa$ odd).
\end{enumerate}
}
\end{rem}
Another type of weighted power variations is given by the
following definition: for $t\in[0,1]$, $f:\R\rightarrow\R$ and
$\kappa\geq 2$, let
\begin{eqnarray*}
S_n^{(\kappa)}(f,t)&=&
\sum_{k=0}^{
\left\lfloor \frac12\big(2^{\frac{n}2}t-1\big)\right\rfloor
} \!\!\!\!f\big(Z(T_{2k+1,n})\big)
\left[\big(Z(T_{2k+2,n})-Z(T_{2k+1,n})\big)^{\kappa}\right.\\
&&\hskip5cm \left.+(-1)^{\kappa+1}
\big(Z(T_{2k+1,n})-Z(T_{2k,n})\big)^{\kappa}\right].
\end{eqnarray*}
This type of variations have been introduced very recently by
Swanson \cite{Swanson} (see, more precisely, relations (1.6) and
(1.7) in \cite{Swanson}), and used in order to obtain a change of
variables formula (in law) for the solution of the stochastic heat
equation driven by a space/time white noise. Since our approach allows us to treat this type of {\sl
signed} weighted power variations, we will also prove the
following result:
\begin{thm}\label{thm2}
Let $f:\R\rightarrow\R$ belong to ${\rm C}^2$ with $f'$ and $f''$ bounded, and
$\kappa\geq 2$ be an integer.
Then, as $n\rightarrow\infty$,
\begin{enumerate}
\item if $\kappa$ is even,
$\left\{X_x,2^{(\kappa-1)\frac{n}4}\,S_n^{(\kappa)}(f,t)\right\}_{x\in\R,\,t\in[0,1]}$
converges in the sense of f.d.d. to
\begin{eqnarray}
\left\{X_x,
\sqrt{\mu_{2\kappa}-\mu_{\kappa}^2}\int_0^{Y_t}f(X_z)dB_z\right\}_{x\in\R,\,t\in[0,1]};
\label{gen2bis}
\end{eqnarray}
\item if $\kappa$ is odd,
$\left\{X_x,
2^{(\kappa-1)\frac{n}4}\,S_n^{(\kappa)}(f,t)\right\}_{x\in\R,\,t\in[0,1]}$
converges in the sense of f.d.d. to
\begin{eqnarray}
\left\{X_x,\int_0^{Y_t} f(X_z)\big(\mu_{\kappa+1}\,
d^\circ X_z+\sqrt{\mu_{2\kappa}-\mu_{\kappa+1}^2}\,dB_z\big)\right\}
_{x\in\R,\,t\in[0,1]}.\label{gen3bis}
\end{eqnarray}
\end{enumerate}
\end{thm}
\begin{remark}
{\rm
\begin{enumerate}
\item See also Burdzy \cite{bur} for an alternative study of (non-weighted) signed variations of I.B.M. Note,
however, that the approach of \cite{bur} is not based on the use of the I.S.S., but rather on thinning
deterministic partitions of the time axis.
\item
The limits and the rates of convergence in
(\ref{gen3}) and (\ref{gen3bis}) are the same, while the limits
and the rates of convergences in (\ref{gen2}) and (\ref{gen2bis})
are different.
\end{enumerate}
}
\end{remark}

The rest of the paper is organized as follows. In Section 2, we
state and prove some ancillary results involving weighted sums of
polynomial transformations of Brownian increments. Section 3 is
devoted to the proof of Theorem \ref{thm1}, while in Section 4 we
deal with Theorem \ref{thm2}.

\section{Preliminaries}

In order to prove Theorem \ref{thm1} and Theorem \ref{thm2}, we
shall need several asymptotic results, involving quantities that
are solely related to the Brownian motion $X$. The aim of this
section is to state and prove these results, that are of clear
independent interest.

We let the notation of the Introduction prevail: in particular,
$X$ and $B$ are two independent two-sided Brownian motions, and
$Y$ is a one-sided Brownian motion independent of $X$ and $B$. For
every $n\geq 1$, we also define the process $X^{(n)}=\{X^{(n)}_t\}_{t\geq
0}$ as
\begin{equation}
X_{t}^{\left( n\right) }=2^{n/4}X_{t2^{-n/2}}.  \label{cosa}
\end{equation}%

\begin{rem}
{\rm In what follows, we will work with the dyadic partition of
order $n/2$, instead of that of order $n$, since the former
emerges very naturally in the proofs of Theorem \ref{thm1} and
Theorem \ref{thm2}, as given, respectively, in Section \ref{S :
TH1.2} and Section \ref{S : TH1.4} below. Plainly, the results
stated and proved in this section can be easily reformulated in
terms of \textsl{any} sequence of partitions with equidistant
points and with meshes converging to zero. }
\end{rem}

The following result plays an important role in this section. In
the next statement, and for the rest of the paper, we will freely
use the language of Wiener chaos and Hermite polynomials. The
reader is referred e.g. to Chapter 1 in \cite{Nbook} for any
unexplained definition or result.
\begin{thm}\label{thm-PT} (Peccati and Tudor \cite{PT}).
Fix $d\geq 2$, fix $d$ natural numbers $1\leq n_1\leq \ldots \leq
n_d$ and, for every $k\geq 1$, let ${\bf
F}^k=(F_1^k,\ldots,F_d^k)$ be a vector of $d$ random variables
such that, for every $j=1,\ldots,d$, the sequence of $F^k_j$,
$k\geq 1$, belongs to the $n_j$th Wiener chaos associated with
$X$. Suppose that, for every $1\leq i,j\leq d$,
$\lim_{k\rightarrow\infty} E(F_i^kF_j^k)=\delta_{ij}$, where
$\delta_{ij}$ is Kronecker symbol. Then, the following two
conditions are equivalent:
\begin{enumerate}
\item[(i)] The sequence ${\bf F}^k$, $k\geq 1$, converges in distribution to a standard centered Gaussian
vector $\mathscr{N}(0,{\bf I}_d)$ (${\bf I}_d$ is the $d\times d$ identity matrix),
\item[(ii)] For every $j=1,\ldots,d$, the sequence $F_j^k$, $k\geq 1$, converges in distribution to a
standard Gaussian random variable $\mathscr{N}(0,1)$.
\end{enumerate}
\end{thm}

The forthcoming proposition is the key to the main results of this
section.

\bigskip

Given a polynomial $\textsf{P}:\mathbb{R}\rightarrow \mathbb{R}$,
we say that $\textsf{P}$ has \textsl{centered Hermite rank} $\geq
2$ whenever $E\left[G\,\textsf{P}(G)\right]=0$ (where $G$ is a standard Gaussian random variable).
Note that $\mathsf{P}$ has centered Hermite rank $\geq 2$ if, and
only if, $\textsf{P}\left( x\right) -E\left[ \textsf{P}\left(
G\right) \right] $
 has Hermite rank $\geq 2$, in the sense of Taqqu \cite{T1}.

\begin{proposition}
\label{P : Incr}Let $\mathsf{P}:\mathbb{R}\rightarrow \mathbb{R}$
be a polynomial with centered Hermite rank $\geq 2$. Let $\alpha
,\beta \in \mathbb{R}$
and denote by $\phi :\mathbb{N}%
\rightarrow \mathbb{R}$ the function defined by the relation:
$\phi (i)$ equals $\alpha $ or $\beta $, according as $i$ is even
or odd. Fix an integer $N\geq 1$ and, for every $j=1,...,N$, set
\begin{eqnarray*}
M_{j}^{\left( n\right) }&=&2^{-n/4}\sum_{i=\lfloor \left( j-1\right)
2^{n/2}\rfloor +1}^{\lfloor j2^{n/2}\rfloor }\phi \left( i\right) \left\{\mathsf{P}\left(
X_{i}^{\left( n\right) }-X_{i-1}^{\left( n\right) }\right)
-E[\mathsf{P}(G)]\right\}
  \text{,}\\
M_{j+N}^{\left( n\right) }&=&2^{-n/4}\sum_{i=\lfloor \left( j-1\right) 2^{n/2}\rfloor +1}
^{\lfloor j2^{n/2}\rfloor }
(-1)^i \left(X_{i}^{(n)}-X_{i-1}^{(n)}\right)
\text{,}
\end{eqnarray*}%
where $G\sim\mathscr{N}(0,1)$ is a standard Gaussian random
variable. Then, as  $n\rightarrow \infty $, the random vector
\begin{equation}
\left\{M_{1}^{\left( n\right) },...,M_{2N}^{(n)} ;
\{X_t\}_{t\geq 0}\right\}   \label{gg}
\end{equation}%
converges weakly in the space $\R^{2N}\times {\rm C}^0(\R_+)$ to
\begin{equation}
\left\{ \sqrt{\frac{\alpha ^{2}+\beta ^{2}}{2}\,\mathbf{Var}
\big(\mathsf{P}(G)\big) }\big( \Delta B\left( 1\right) ,...,\Delta
B\left( N\right) \big) ; \Delta B\left( N+1\right) ,...,\Delta
B\left( 2N\right); \{X_t\}_{t\geq 0}\right\}   \label{g}
\end{equation}%
where $\Delta B\left( i\right) =B_i -B_{i-1}$, $i=1,...,N$.
\end{proposition}

\begin{proof}

For the sake of notational simplicity, we provide the proof only in the case
where $\alpha =\sqrt{2}$ and $\beta =0$, the extension to the general case
being a straightfroward consequence of the independence of the Brownian
increments. For every $h\in L^{2}\left( \R_+\right) $, we
write $X\left( h\right) =\int_{0}^{\infty}h\left( s\right) dX_{s}$. To prove the
result it is sufficient to show that, for every $\mathbf{\lambda }=\left(
\lambda _{1},...,\lambda _{2N+1}\right)
\in \mathbb{R}^{2N+1}$
and every $h\in L^{2}\left(\R_+\right) $, the sequence of random variables
\begin{equation*}
F_{n} =
\sum_{j=1}^{2N}%
\lambda _{j}M_{j}^{\left( n\right) }
+\lambda _{2N+1}X\left( h\right)
\end{equation*}%
converges in law to
\begin{equation*}
\sqrt{\mathbf{Var}\big(\mathsf{P}(G)\big) }
\sum_{j=1}^{N}\lambda_{j}\Delta B\left( j\right)
+\sum_{j=N+1}^{2N}\lambda_{j}\Delta B\left(j\right) +\lambda
_{2N+1}X\left( h\right) \text{.}
\end{equation*}%
We start by observing that the fact that $\mathsf{P}$ has centered Hermite rank $%
\geq 2$
implies that $%
\mathsf{P}$ is such that
\begin{equation}
\mathsf{P}\left( X_{i}^{\left( n\right) }-X_{i-1}^{\left( n\right)
}\right) -E\left[ \mathsf{P}\left( G\right) \right]
 =\sum_{m=2}^{\kappa
}b_{m}H_{m}\left( X_{i}^{\left( n\right) }-X_{i-1}^{\left( n\right) }\right)
\text{,}\quad\mbox{for some $\kappa\geq 2$},  \label{cm}
\end{equation}%
where $H_{m}$ denotes the $m$th Hermite polynomial, and the
coefficients $b_{m}$ are real-valued and uniquely determined by
(\ref{cm}). Moreover, one has that
\begin{equation}
\mathbf{Var}\big(\mathsf{P}(G)\big) =\sum_{m=2}^{\kappa }b_{m}^{2}\,\mathbb{%
E}\left[ H_{m}\left(G\right) ^{2}\right] =\sum_{m=2}^{\kappa
}b_{m}^{2}\,m!\text{.}  \label{var}
\end{equation}%
We can now write
\begin{eqnarray*}
F_{n} &=&\lambda _{2N+1}X\left( h\right) +%
\sqrt{2}\,\,\,2^{-\frac{n}{4}}\sum_{j=1}^{N}\lambda _{j}\sum_{\substack{ %
i=\lfloor \left( j-1\right) 2^{n/2}\rfloor +1 \\ i\,\,\text{even}}}^{\lfloor j2^{n/2}\rfloor }\,\,\sum_{m=2}^{%
\kappa }b_{m}H_{m}\left( X_{i}^{\left( n\right) }-X_{i-1}^{\left( n\right)
}\right) \\
&&\hskip5cm+2^{-\frac{n}4}\sum_{j=1}^{N}\lambda_{N+j}
\sum_{i=\lfloor \left( j-1\right) 2^{n/2}\rfloor +1}^{\lfloor j2^{n/2}\rfloor }
(-1)^i \left(X^{(n)}_{i}-X^{(n)}_{i-1}\right)
 \\
&=&\lambda _{2N+1}X\left( h\right) +\sum_{m=2}^{\kappa
}b_{m}\sum_{j=1}^{N}\lambda _{j}\,\,\sqrt{2}\,\,\,2^{-\frac{n}{4}}\sum
_{\substack{ i=\lfloor \left( j-1\right) 2^{n/2}\rfloor +1 \\ i\,\,\text{even}}}%
^{\lfloor j2^{n/2}\rfloor }\,\,H_{m}\left( X_{i}^{\left( n\right) }-X_{i-1}^{\left( n\right)
}\right) \\
&&\hskip5cm+2^{-\frac{n}4}\sum_{j=1}^{N}\lambda_{N+j}
\sum_{i=\lfloor \left( j-1\right) 2^{n/2}\rfloor +1}^{\lfloor j2^{n/2}\rfloor }
(-1)^i \left(X^{(n)}_{i}-X^{(n)}_{i-1}\right).
\end{eqnarray*}%
By using the independence of the Brownian increments, the Central Limit
Theorem and Theorem \ref{thm-PT}, we deduce that the $%
\kappa+1 $ dimensional vector
\begin{eqnarray*}
&&\left\{ X\left( h\right) ;
2^{-\frac{n}4}\sum_{j=1}^{N}\lambda_{N+j}
\sum_{i=\lfloor \left( j-1\right) 2^{n/2}\rfloor +1}^{\lfloor j2^{n/2}\rfloor }
(-1)^i \left(X^{(n)}_{i}-X^{(n)}_{i-1}\right);\right.\\
&&\hskip5cm\left.
\sum_{j=1}^{N}\lambda _{j}\,\,\sqrt{2}\,\,\,2^{-%
\frac{n}{4}}\sum_{\substack{ i=\lfloor \left( j-1\right) 2^{n/2}\rfloor +1 \\
i\,\,\text{even}
}}^{\lfloor j2^{n/2}\rfloor }\,\,H_{m}\left( X_{i}^{\left( n\right) }-X_{i-1}^{\left( n\right)
}\right) :m=2,...,\kappa
\right\}
\end{eqnarray*}%
converges in law to
\begin{equation*}
\left\{ \left\Vert h\right\Vert_2 \times G_{0};\sum_{j=1}^{N}\lambda _{N+j}\,
G_{j,1};\sum_{j=1}^{N}\lambda _{j}\,\sqrt{%
m!}G_{j,m}
:m=2,...,\kappa \right\} \text{,}
\end{equation*}%
where $\left\{ G_{0};G_{j,m}:j=1,...,N\text{, \ }m=1,...,\kappa \right\} $
is a collection of i.i.d. standard Gaussian random variables $\mathscr{N}(0,1)$.
This implies that $F_{n} $ converges in
law, as $n\rightarrow \infty $, to%
\begin{eqnarray*}
&&\lambda _{2N+1}\left\Vert h\right\Vert_2 G_{0}
+\sum_{j=1}^{N}\lambda_{N+j}G_{j,1}
+\sum_{j=1}^{N}\lambda
_{j}\sum_{m=2}^{\kappa }b_{m}\sqrt{m!}G_{j,m}\\
&\overset{\mathrm{Law}}{=}& \lambda _{2N+1}X\left( h\right)
+\sum_{j=N+1}^{2N}\lambda_{j}\Delta B\left( j\right)
+\sum_{j=1}^{N}\lambda _{j}\sqrt{\mathbf{Var}%
\big(\mathsf{P}(G)\big) }\Delta B\left( j\right) \text{,}
\end{eqnarray*}%
where we have used (\ref{var}). This proves our claim.
\end{proof}

\begin{rem}
\rm{ It is immediately verified that the sequence of Brownian
motions appearing in (\ref{cosa}) is asymptotically independent of
$X$ (just compute the covariance function of the 2-dimensional
Gaussian process $(X,X^{(n)})$). However, by inspection of the
proof of Proposition \ref{P : Incr}, one sees that this fact is
immaterial in the proof of the asymptotic independence of the
vector $(M_1^{n},...,M_{2N}^{n})$ and $X$. Indeed, such a result
\textsl{depends uniquely of the fact that the polynomial
$\mathsf{P}-E(\mathsf{P}(G))$ has an Hermite rank of order
strictly greater than one}. This is a consequence of Theorem
\ref{thm-PT}. It follows that the statement of Proposition \ref{P
: Incr} still holds when the sequence $\{X^{(n)}\}$ is replaced by
a sequence of Brownian motions $\{X^{(*,n)}\}$ of the type
$$X^{(*,n)}(t)=\int\psi_n(t,z)dX_z, \quad t\geq 0, \quad n \geq 1,
$$ where, for each $t$, $\psi_n$ is a square-integrable
deterministic kernel.

}
\end{rem}

Once again, we stress that $d$ and $d^\circ$ denote, resp., the Wiener-It\^o integral
and the Stratonovich integral, see also Remark \ref{rk13} (1).

\begin{theorem}
\label{T : BastaJacod}Let the notation and assumptions of
Proposition \ref{P : Incr} prevail (in particular, $\mathsf{P}$
has centered Hermite rank $\geq 2$), and set
\begin{eqnarray*}
J_{t}^{\left( n\right) }\left( f\right) &=&2^{-\frac{n}{4}}\,\frac12\sum_{j=1}^{\left%
\lfloor 2^{n/2}t\right\rfloor }
\big(f(X_{\left( j-1\right)2^{-n/2}})+f(X_{j2^{-n/2}})\big)
\\
&&\hskip10mm\times\left[\phi \left( j\right) \left\{ \mathsf{P}\left(
X_{j}^{\left( n\right)
}-X_{j-1}^{\left( n\right) }\right) -E\left[ \mathsf{P}\left( G\right) \right] \right\}
+\gamma(-1)^j\left( X_{j}^{\left( n\right)
}-X_{j-1}^{\left( n\right) }\right)
 \right]
\text{, \ }t\geq 0 \text{,}
\end{eqnarray*}%
where $\gamma\in\R$ and the real-valued function $f$ is supposed to belong to ${\rm C}^2$
with $f'$ and $f''$ bounded. Then, as $%
n\rightarrow +\infty $, the two-dimensional process%
\begin{equation*}
\left\{ J_{t}^{\left( n\right) }\left( f\right) ,
X_{t}\right\}_{t\geq 0}
\end{equation*}%
converges in the sense of f.d.d. to
\begin{equation}\label{z1}
\left\{ \sqrt{\gamma^2+\frac{\alpha ^{2}+\beta
^{2}}{2}\mathbf{Var}\big(\mathsf{P}(G)\big)}\, \int_{0}^{t}f\left(
X_{s}\right) dB_s, X_{t}\right\}_{t\geq 0} \text{.}
\end{equation}
\end{theorem}

\begin{proof}

Set $\sigma :=\sqrt{\frac{\alpha ^{2}+\beta
^{2}}{2}\mathbf{Var}\big(\mathsf{P}(G)\big) }$. 
We have
$$
J_t^{(m)}(f)=\widetilde{J}_t^{(m)}(f)+r_t^{(m)}(f)+s_t^{(m)}(f),
$$
where
\begin{eqnarray*}
\widetilde{J}_t^{(m)}(f)&=&2^{\frac{m}4}\sum_{j=1}^{\lfloor 2^{\frac{m}2}t\rfloor}
f(X_{(j-1)2^{-\frac{m}2}})\left[ \phi \left( j\right) \left\{
\mathsf{P}\left( X_{j}^{\left( m\right) }-X_{j-1}^{\left( m\right)
}\right) -E\left[ \mathsf{P}\left( G\right) \right] \right\}\right.\\
&&\hskip7cm\left.
+\gamma(-1)^j\left( X_{j}^{\left( m\right)
}-X_{j-1}^{\left( m\right) }\right)
 \right]\\
r_t^{(m)}(f)&=&2^{-\frac{m}2}\,\frac12\sum_{j=1}^{\lfloor 2^{\frac{m}2}t\rfloor}
f'(X_{(j-1)2^{-\frac{m}2}})\left( X_{j}^{\left( m\right)
}-X_{j-1}^{\left( m\right) }\right)\\
&&\hskip1cm\times
\left[ \phi \left( j\right) \left\{
\mathsf{P}\left( X_{j}^{\left( m\right) }-X_{j-1}^{\left( m\right)
}\right) -E\left[ \mathsf{P}\left( G\right) \right] \right\}
+\gamma(-1)^j\left( X_{j}^{\left( m\right)
}-X_{j-1}^{\left( m\right) }\right)
 \right]\\
s_t^{(m)}(f)&=&2^{-\frac{3m}4}\,\frac14\sum_{j=1}^{\lfloor 2^{\frac{m}2}t\rfloor}
f''(X_{\theta_{j,m}})\left( X_{j}^{\left( m\right)
}-X_{j-1}^{\left( m\right) }\right)^2\\
&&\hskip1cm\times
\left[ \phi \left( j\right) \left\{
\mathsf{P}\left( X_{j}^{\left( m\right) }-X_{j-1}^{\left( m\right)
}\right) -E\left[ \mathsf{P}\left( G\right) \right] \right\}
+\gamma(-1)^j\left( X_{j}^{\left( m\right)
}-X_{j-1}^{\left( m\right) }\right)
 \right],
\end{eqnarray*}
for some $\theta_{j,m}$ between $(j-1)2^{-\frac{m}2}$ and $j2^{-\frac{m}2}$. We decompose
\begin{eqnarray*}
r_t^{(m)}(f)&=&2^{-\frac{m}2}\,\frac12\sum_{j=1}^{\lfloor 2^{\frac{m}2}t\rfloor}
\phi \left( j\right)
f'(X_{(j-1)2^{-\frac{m}2}})\left( X_{j}^{\left( m\right)
}-X_{j-1}^{\left( m\right) }\right)
\left\{
\mathsf{P}\left( X_{j}^{\left( m\right) }-X_{j-1}^{\left( m\right)
}\right) -E\left[ \mathsf{P}\left( G\right) \right] \right\}\\
&&+\frac{\gamma}2\,2^{-\frac{m}2}\sum_{j=1}^{\lfloor 2^{\frac{m}2}t\rfloor}
(-1)^j
f'(X_{(j-1)2^{-\frac{m}2}})
\left( X_{j}^{\left( m\right)
}-X_{j-1}^{\left( m\right) }\right)
^2=r_t^{(1,m)}(f)+r_t^{(2,m)}(f).
\end{eqnarray*}
By independence of increments and because $E\left[G\big(\mathsf{P}(G)-E[\mathsf{P}(G)]\big)\right]=0$
we have
$$
E\big|r_t^{(1,m)}(f)\big|^2=\frac14\,2^{-m}\,E\left[G^2\big(\mathsf{P}(G)-E[\mathsf{P}(G)]\big)^2\right]
\sum_{j=1}^{\lfloor 2^{\frac{m}2}t\rfloor} \phi^2(j) E\big|f'(X_{(j-1)2^{-\frac{m}2}})\big|^2 = O(2^{-\frac{m}2}).
$$
For $r_t^{(2,m)}(f)$, we can decompose
\begin{eqnarray*}
r_t^{(2,m)}(f)&=&\frac{\gamma}2\,2^{-\frac{m}2}
\sum_{j=1}^{\lfloor 2^{\frac{m}2}t\rfloor} (-1)^j f'(X_{(j-1)2^{-\frac{m}2}})\left[
\left( X_{j}^{\left( m\right)
}-X_{j-1}^{\left( m\right) }\right)^2-1
\right]\\
&&+\frac{\gamma}2\,2^{-\frac{m}2}
\sum_{j=1}^{\lfloor 2^{\frac{m}2}t\rfloor} (-1)^j f'(X_{(j-1)2^{-\frac{m}2}})
=r_t^{(2,1,m)}(f)+r_t^{(2,2,m)}(f).
\end{eqnarray*}
By independence of increments and because $E\big(G^2-1\big)=0$, we have
$$
E\big|r_t^{(2,1,m)}(f)\big|^2=\frac{\gamma^2}4\,\mathbf{Var}(G^2)\,2^{-m}\,
\sum_{j=1}^{\lfloor 2^{\frac{m}2}t\rfloor} E\big|f'(X_{(j-1)2^{-\frac{m}2}})\big|^2 = O(2^{-\frac{m}2}).
$$
For $r_t^{(2,2,m)}(f)$, remark that
$$
2^{-\frac{m}2} \sum_{k=0}^{\lfloor 2^{\frac{m}2}t\rfloor}
(-1)^j f'(X_{(j-1)2^{-\frac{m}2}}) = 2^{-\frac{m}2} \sum_{k=0}^{\lfloor \lfloor 2^{\frac{m}2}t\rfloor/2\rfloor}
\big( 
f'(X_{(2k+1)2^{-\frac{m}2}})
-f'(X_{(2k)2^{-\frac{m}2}})  
\big)
+O(2^{-\frac{m}2})
$$
so that $E\big|r_t^{(2,2,m)}(f)\big|=O(2^{-\frac{m}4})$. Thus, we obtained that
$E\big| r_t^{(m)}(f)\big|\longrightarrow 0$ as $m\to\infty$.
Since we obviously have, using the boundedness of $f''$, that
$$
E\big| s_t^{(m)}(f)\big|\longrightarrow 0\quad\mbox{as $m\to\infty$},
$$
we see that the convergence result in Theorem \ref{T : BastaJacod} is equivalent to the
convergence of the pair $\left\{ \widetilde{J}_t^{(n)}(f),X_t\right\}_{t\geq 0}$ to the object in (\ref{z1}).

Now, for every $m\geq n$, one
has that%
\begin{eqnarray*}
J_{t}^{\left( m\right) }\left( f\right)  &=&2^{-\frac{m}{4}%
}\sum_{j=1}^{\left\lfloor 2^{n/2}t\right\rfloor
}\!\!\sum_{i=\lfloor \left( j-1\right)
2^{\frac{m-n}{2}}\rfloor +1}^{\lfloor j2^{\frac{m-n}{2}}\rfloor } \!\!\!\!\!\!
f(X_{\left(i-1\right) 2^{-m/2}})
\\
&&\hskip1.5cm\times \left[ \phi \left( i\right) \left\{
\mathsf{P}\left( X_{i}^{\left( m\right) }-X_{i-1}^{\left( m\right)
}\right) -E\left[ \mathsf{P}\left( G\right) \right] \right\}
+\gamma(-1)^i\left( X_{i}^{\left( m\right)
}-X_{i-1}^{\left( m\right) }\right)
 \right]
\\
&=&A_{t}^{\left( m,n\right) }+B_{t}^{\left( m,n\right) },
\end{eqnarray*}%
where
\begin{eqnarray*}
A_{t}^{\left( m,n\right) } &=&2^{-\frac{m}{4}}\sum_{j=1}^{\left\lfloor
2^{n/2}t\right\rfloor }f(X_{\left( j-1\right) 2^{-n/2}})\\
&&\hskip5mm\times\sum_{i=\lfloor \left(
j-1\right) 2^{\frac{m-n}2}\rfloor +1}^{\lfloor j2^{\frac{m-n}2}\rfloor }\!\!\!\left[\phi \left( i\right)
\left\{ \mathsf{P}\left( X_{i}^{\left( m\right) }-X_{i-1}^{\left( m\right) }\right) -\mathbb{E%
}\left[ \mathsf{P}\left( G\right) \right] \right\}
+\gamma(-1)^i\left( X_{i}^{\left( m\right)
}-X_{i-1}^{\left( m\right) }\right)
 \right]\\
B_{t}^{\left( m,n\right) }
&=&2^{-\frac{m}{4}}\sum_{j=1}^{\left\lfloor
2^{n/2}t\right\rfloor }\sum_{i=\lfloor \left( j-1\right)
2^{\frac{m-n}2}\rfloor +1} ^{\lfloor j2^{\frac{m-n}2}\rfloor
} \left[ f(X_{\left(
i-1\right) 2^{-m/2}})-f(X_{\left( j-1\right) 2^{-n/2}})\right]
  \\
&& \text{\ \ \ \ \ \ \ \ \ \ \ \ \ \ \ \ } \times \left[\phi
\left( i\right) \left\{ \mathsf{P}\left( X_{i}^{\left( m\right)
}-X_{i-1}^{\left( m\right) }\right) -E\left[ \mathsf{P}\left(
G\right) \right] \right\} +\gamma(-1)^i\left( X_{i}^{\left(
m\right) }-X_{i-1}^{\left( m\right) }\right) \right].
\end{eqnarray*}%
We shall study $A^{\left( m,n\right) }$ and $B^{\left( m,n\right) }$
separately. By Proposition \ref{P : Incr}, we know that, as $m\rightarrow
\infty $, the random element%
\begin{eqnarray*}
&&\left\{ X;2^{-\frac{m}{4}}\sum_{i=\lfloor \left( j-1\right)
2^{\frac{m-n}2}\rfloor +1}^{\lfloor j2^{\frac{m-n}2}\rfloor }\phi \left( i\right) \left\{
\mathsf{P}\left( X_{i}^{\left( m\right) }-X_{i-1}^{\left( m\right)
}\right) -E\left[ \mathsf{P}\left( G\right) \right] \right\}
:j=1,...,2^{n/2}
;\right.\\
&&\hskip4cm\left.
2^{-\frac{m}{4}}\sum_{i=\lfloor \left( j-1\right)
2^{\frac{m-n}2}\rfloor +1}^{\lfloor j2^{\frac{m-n}2}\rfloor }(-1)^i \left( X_{i}^{\left( m\right) }-
X_{i-1}^{\left( m\right) }\right)
:j=1,...,2^{n/2}
\right\}
\end{eqnarray*}%
converges in law to
\begin{eqnarray*}
&&\left\{ X;
2^{-\frac{n}{4}}
\sigma \,\Delta
B\left( j\right) :j=1,...,2^{n/2};
2^{-\frac{n}{4}}\,\Delta
B\left( j+2^{n/2}\right) :j=1,...,2^{n/2}
\right\}  \\
&\overset{\mathrm{Law}}{=}&\left\{ X;\sigma %
\left( B\left( j2^{-n/2}\right) -B\left( \left( j-1\right)
2^{-n/2}\right) \right) :j=1,...,2^{n/2};\right.\\
&&\hskip6cm \left.B_{2}\left( j2^{-n/2}\right) -B_{2}\left( \left( j-1\right)
2^{-n/2}\right) :j=1,...,2^{n/2}
\right\},
\end{eqnarray*}%
where $B_2$ denotes a standard Brownian motion, independent of $X$ and $B$.
Hence, as $m\rightarrow \infty $,
$$
\left\{ X;A^{\left( m,n\right) }\right\}
\,\,\,\overset{{\rm f.d.d.}}{\Longrightarrow}\,\,\,
\left\{ X;A^{\left( \infty,n\right) }\right\}
$$
where
\begin{eqnarray*}
&&A_{t}^{(\infty,n)}:=
\sigma\sum_{j=1}^{\left\lfloor 2^{n/2}t\right\rfloor }
f(X_{\left( j-1\right) 2^{-n/2}})%
\left[ B\left( j2^{-n/2}\right) -B\left( \left( j-1\right)
2^{-n/2}\right) \right] \\
&&\hskip2cm
+\gamma\sum_{j=1}^{\left\lfloor 2^{n/2}t\right\rfloor }
f(X_{\left( j-1\right) 2^{-n/2}})%
\left[ B_{2}\left( j2^{-n/2}\right) -B_{2}\left( \left( j-1\right)
2^{-n/2}\right) \right].
\end{eqnarray*}%
By letting $n\rightarrow \infty $, and by using the independence
of $X, B$ and $B_2$, one obtains that $A^{(\infty,n)}$ converges
uniformly on compacts in probability (u.c.p.) towards
$$
A_{t}^{(\infty,\infty)}\triangleq \int_{0}^{t}f\left( X_{s}\right)
\big(\sigma\,dB(s)+\gamma\,dB_2(s)\big).
$$
This proves that, by letting $m$ and then $n$ go to
infinity $\left\{ X;A^{\left( m,n\right) }\right\} $ converges
in the sense of f.d.d. to
$$
\left\{ X;A^{(\infty,\infty)}\right\} \overset{\mathrm{Law}}{=}
\left\{X; \sqrt{\sigma^2 +\gamma^2}\, \int_{0}^{\cdot}f\left(
X_{s}\right)dB(s)\right\} .$$

To conclude the proof of the Theorem we shall show that, by
letting $m$ and then $n$ go to infinity,  $ B_{t}^{\left(
m,n\right) }$ converges to zero in $L^{2}$, for any fixed $t>0$.
To see this, write $K=\max(|\alpha|,|\beta|)$ and observe that,
for every $t>0$, the independence of the Brownian increments 
yields the chain of inequalities:
\begin{eqnarray*}
E\left[ \left\vert B_{t}^{\left( m,n\right) }\right\vert^2 \right]
&=&2^{-m/2}\sum_{j=1}^{ \left\lfloor 2^{n/2}t\right\rfloor
}\sum_{i=\lfloor \left( j-1\right) 2^{\frac{m-n}2}\rfloor +1}^{\lfloor j2^{\frac{m-n}2}\rfloor }
E\left\{\left|
f(X_{\left( i-1\right) 2^{-m/2}})-f(X_{\left( j-1\right) 2^{-n/2}})
\right|^2\right\}\\
&&\times E\left\{\left|\phi \left( i\right) \left\{
\mathsf{P}\left( X_{i}^{\left( m\right) }-X_{i-1}^{\left( m\right)
}\right) -E\left[ \mathsf{P}\left( G\right) \right] \right\}
+\gamma(-1)^i\left( X_{i}^{\left( m\right)
}-X_{i-1}^{\left( m\right) }\right)\right|^2\right\}\\
&\leq & 2|f'|^2_\infty\big(K^2\mathbf{Var} \big(\mathsf{P}(G)\big)
+\gamma^2\big) 2^{-m/2}\!\!\sum_{j=1}^{\left\lfloor
2^{n/2}t\right\rfloor }\!\!\! \sum_{i=\lfloor \left( j-1\right)
2^{\frac{m-n}2}\rfloor +1}
^{\lfloor j2^{\frac{m-n}2}\rfloor }%
\!\!\!\!\!\!\left[ \left( i-1\right) 2^{-\frac{m}2}\!-\!\left( j-1\right) 2^{-\frac{n}2}\right]  \\
&\leq &2|f'|^2_\infty\big(K^2\mathbf{Var} \big(\mathsf{P}(G)\big)
+\gamma^2\big) 2^{-\frac{m+n}2 }\sum_{j=1}^{\left\lfloor
2^{n/2}t\right\rfloor }\sum_{i=\lfloor \left(
j-1\right) 2^{\frac{m-n}2}\rfloor +1}^{\lfloor j2^{\frac{m-n}2}\rfloor }1\\
&\leq&2|f'|^2_\infty t\big(K^2\mathbf{Var} \big(\mathsf{P}(G)\big)
+\gamma^2\big)2^{-n/2}.
\end{eqnarray*}%
This shows that
\begin{equation*}
\limsup_{m\rightarrow \infty }E\left[\left\vert B_{t}^{\left(
m,n\right) }\right\vert^2 \right]\leq {\rm cst.}\,2^{-n/2},
\end{equation*}%
and the desired conclusion is obtained by letting $n\rightarrow \infty .$
\end{proof}

\bigskip

We now state several consequences of Theorem \ref{T : BastaJacod}.
The first one (see also Jacod \cite{Jacod}) is obtained by setting
$\alpha=\beta=1$ (that is, $\phi$ is identically one), $\gamma=0$
and by recalling that, for an even integer $\kappa \geq 2$, the
polynomial $P(x)=x^{\kappa }-\mu_\kappa$ has centered Hermite rank
$\geq 2$.

\bigskip

\begin{corollary}\label{cor9}
Let $f:\mathbb{R}\rightarrow\mathbb{R}$ belong to ${\rm C}^2$ with $f'$ and $f''$ bounded, and fix an even
integer $\kappa \geq 2 $. For $t\geq 0$, we set:
\begin{equation}
J_{t}^{\left( n\right) }\left( f\right) =2^{-\frac{n}{4}}\,\frac12\sum_{j=1}^{\left%
\lfloor 2^{n/2}t\right\rfloor }
\big(f(X_{\left( j-1\right) 2^{-n/2}})+f(X_{j2^{-n/2}})
\big)\left\{
2^{\kappa \frac{n}{4}}\left( X_{j2^{-n/2}}-X_{\left( j-1\right) 2^{-n/2}}\right)
^{\kappa }-\mu _{\kappa }\right\}.
\label{Jf1}
\end{equation}%
Then, as $n\rightarrow +\infty $, $\left\{ J_{t}^{\left( n\right)
}\left( f\right) ,X_{t}\right\}_{t\geq 0}$ converges in the sense
of f.d.d. to
\begin{equation*}
\left\{ \sqrt{\mu _{2\kappa }-\mu _{\kappa }^{2}}\int_{0}^{t}f\left(
X_{s}\right) dB_s ,X_{t}\right\}_{t\geq 0}
\text{.}
\end{equation*}
\end{corollary}

The next result derives from Theorem \ref{T : BastaJacod} in the case $%
\alpha=1$, $\beta=-1$ and $\gamma=0$.

\bigskip

\begin{corollary}\label{cor10}
Let $f:\mathbb{R}\rightarrow\mathbb{R}$ belong to ${\rm C}^2$ with $f'$ and $f''$ bounded,
$\kappa \geq 2$ be an even
integer, and set, for $t\geq 0$:
\begin{equation}
J_{t}^{\left( n\right) }\left( f\right) =2^{(\kappa-1)\frac{n}{4}}\,\frac12\sum_{j=1}^{\left%
\lfloor 2^{n/2}t\right\rfloor }
\big(f(X_{\left( j-1\right) 2^{-n/2}})+f(X_{j2^{-n/2}})\big)
\left( -1\right)
^{j}\left( X_{j2^{-n/2}}-X_{\left( j-1\right)
2^{-n/2}}\right) ^{\kappa }.
\label{Jf-1}
\end{equation}%
Then, as $n\rightarrow +\infty $, the process $\left\{
J_{t}^{\left( n\right) }\left( f\right) ,X_{t}\right\}_{t\geq 0}$
converges in the sense of f.d.d. to
\begin{equation}
\left\{ \sqrt{\mu _{2\kappa }-\mu _{\kappa }^{2}}\int_{0}^{t}f\left(
X_{s}\right) dB_s ,X_{t} \right\}_{t\geq 0}
\text{.}  \label{ob}
\end{equation}
\end{corollary}

\begin{proof}

It is not difficult to see that the convergence result in the statement is
equivalent to the convergence of the pair $\left\{ Z_{t}^{\left( n\right)
}\left( f\right) ,X_{t}\right\}_{t\geq 0}$ to the object in (%
\ref{ob}), where
\begin{equation*}
Z_{t}^{\left( n\right) }\left( f\right) =2^{-\frac{n}{4}}\,\frac12\sum_{j=1}^{\left%
\lfloor 2^{n/2}t\right\rfloor }
\big(f(X_{\left( j-1\right) 2^{-n/2}})+f(X_{j2^{-n/2}})\big)
\left( -1\right)
^{j}\left\{ 2^{\kappa \frac{n}{4}}\left( X_{j2^{-n/2}}-X_{\left( j-1\right)
2^{-n/2}}\right) ^{\kappa }-\mu _{\kappa }\right\} \text{,}
\end{equation*}%
so that the conclusion is a direct consequence of Theorem \ref{T :
BastaJacod}. Indeed, we have
\begin{eqnarray*}
2^{-\frac{n}{4}}\,\frac12\sum_{j=1}^{\left\lfloor 2^{n/2}t\right\rfloor
}(-1)^{j}\big(f(X_{\left( j-1\right) 2^{-n/2}}) +f(X_{j2^{-n/2}})\big)
&=&2^{-\frac{n}{4}}\,\frac12\big(
f(X_{2\lfloor 2^{n/2}t \rfloor 2^{-n/2}})-f(X_0)
\big)
\end{eqnarray*}%
which tends to zero in ${\rm L}^2$, as $n\rightarrow\infty$.
\end{proof}

\bigskip

A slight modification of Corollary \ref{cor10} yields:

\begin{corollary}\label{cor11}
Let $f:\mathbb{R}\rightarrow\mathbb{R}$ belong to ${\rm C}^2$ with $%
f^{\prime }$ and $f^{\prime \prime }$ bounded, let $\kappa \geq 2$
be an even integer, and set:
\begin{eqnarray}
\widetilde{J}_{t}^{\left( n\right) }\left( f\right) &=&2^{(\kappa-1)\frac{n}{4}}
\!\!\!\!\sum_{j=1}^{\left%
\lfloor \frac12\big(2^{\frac{n}2}t-1\big)\right\rfloor }\!\!\!\!
f(X_{\left( 2j+1\right) 2^{-\frac{n}2}})\left[
\left( X_{(2j+2)2^{-\frac{n}2}}-X_{\left(2j+1\right)
2^{-\frac{n}2}}\right) ^{\kappa } \right.\nonumber\\
&&\hskip6cm\left.
-\left( X_{(2j+1)2^{-\frac{n}2}}-X_{\left( 2j\right)
2^{-\frac{n}2}}\right) ^{\kappa } \right],\quad t\geq 0.\nonumber
\end{eqnarray}%
Then, as $n\rightarrow +\infty $, the process $\left\{
\widetilde{J}_{t}^{\left( n\right) }\left( f\right)
,X_{t}\right\}_{t\geq 0}$ converges in the sense of f.d.d. to
(\ref{ob}).
\end{corollary}

\begin{proof}

By separating the sum according to the eveness of $j$, one can write
$$
\widetilde{J}_{t}^{\left( n\right) }\left( f\right)= J_{t}^{\left(
n\right) }\left( f\right) - r_t^{(n)}(f) + s_t^{(n)}(f),
$$
for $J_{t}^{\left( n\right) }\left( f\right)$ defined by (\ref{Jf-1}) and
\begin{eqnarray*}
r_t^{(n)}(f)&=&\frac{2^{(\kappa-1)\frac{n}4}}{2}\sum_{j=1}^{\lfloor
\frac12\big( 2^{\frac{n}2}t-1 \big)\rfloor}\big(
f(X_{(2j+2)2^{-n/2}})-f(X_{(2j+1)2^{-n/2}}) \big) \big(
X_{(2j+2)2^{-n/2}}-X_{(2j+1)2^{-n/2}}
\big)^\kappa\\
s_t^{(n)}(f)&=&\frac{2^{(\kappa-1)\frac{n}4}}{2}\sum_{j=1}^{\lfloor
\frac12\big( 2^{\frac{n}2}t-1 \big)\rfloor}\big(
f(X_{(2j)2^{-n/2}})-f(X_{(2j+1)2^{-n/2}}) \big) \big(
X_{(2j)2^{-n/2}}-X_{(2j+1)2^{-n/2}} \big)^\kappa.
\end{eqnarray*}
We decompose
\begin{eqnarray*}
r_t^{(n)}(f)&=&\frac{2^{(\kappa-1)\frac{n}4}}{2}\sum_{j=1}^{\lfloor
\frac12\big( 2^{\frac{n}2}t-1 \big)\rfloor} f'(X_{(2j+1)2^{-n/2}})
\big( X_{(2j+2)2^{-n/2}}-X_{(2j+1)2^{-n/2}}
\big)^{\kappa+1}\\
&+&\frac{2^{(\kappa-1)\frac{n}4}}{2}\sum_{j=1}^{\lfloor
\frac12\big( 2^{\frac{n}2}t-1 \big)\rfloor}
\!\!\!\!f''(X_{\theta_{j,n}}) \big(
X_{(2j+2)2^{-n/2}}-X_{(2j+1)2^{-n/2}} \big)^{\kappa+2}=
r_t^{(1,n)}(f) +r_t^{(2,n)}(f),
\end{eqnarray*}
for some $\theta_{j,n}$ between $(2j+1)2^{-n/2}$ and $(2j+2)2^{-n/2}$.
By independence of increments and because $E\big[G^{\kappa+1}\big]=0$, we have
$$
E\big|r_t^{(1,n)}(f)\big|^2 = E\big[G^{2\kappa+2}\big]\,2^{-n}\,\frac14\sum_{j=1}^{\lfloor \frac12\big(
2^{\frac{n}2}t-1
\big)\rfloor}
 E\big|
f'(X_{(2j+1)2^{-n/2}})
\big|^2 = O(2^{-n/2}).
$$
For $r_t^{(2,n)}(f)$, we have
$$
E\big|r_t^{(2,n)}(f)\big| = O(2^{-n/4}).
$$
Similarly, we prove that $E\big|s_t^{(n)}(f)\big|$ tends to zero as $n\rightarrow\infty$,
so that the conclusion is a direct consequence of Corollary \ref{cor10}.
\end{proof}

\bigskip
The subsequent results focus on odd powers.

\bigskip

\begin{corollary}\label{cor-odd}
Let $f:\mathbb{R}\rightarrow\mathbb{R}$ belong to ${\rm C}^2$ with $f'$ and $f''$ bounded, $\kappa \geq 3 $ be
an odd integer, and define $J_{t}^{\left( n\right) }\left(
f\right) $ according to (\ref{Jf1}) (remark however that
$\mu_\kappa=0$). Then, as $n\rightarrow +\infty $, $\left\{
J_{t}^{\left( n\right) }\left( f\right) ,X_{t} \right\}_{t\geq 0}
$ converges in the sense of f.d.d. to
\begin{equation}\label{st}
\left\{ \int_{0}^{t}f\left( X_{s}\right) \big( \mu _{\kappa +1}d^\circ X\left(
s\right) +\sqrt{\mu _{2\kappa }-\mu_{\kappa+1}^2}dB_s \big) ,X_{t} \right\}_{t\geq 0} \text{.}
\end{equation}
\end{corollary}

\begin{proof}

One can write:
\begin{eqnarray*}
J_{t}^{\left( n\right) }\left( f\right) &=&2^{-\frac{n}{4}%
}\,\frac12\sum_{j=1}^{\left\lfloor 2^{n/2}t\right\rfloor }
\big(f(X_{\left( j-1\right)2^{-n/2}})+f(X_{j2^{-n/2}})\big)\\
&& \quad\quad\quad\quad\quad\quad\quad\quad\quad\quad\quad \times
\left\{ \left( X_{j}^{\left( n\right) }-X_{\left( j-1\right)
}^{\left( n\right) }\right) ^{\kappa }-\mu _{\kappa +1}\left(
X_{j}^{\left(
n\right) }-X_{\left( j-1\right) }^{\left( n\right) }\right) \right\} \\
&&+\frac{\mu _{\kappa +1}}{2}\sum_{j=1}^{\left\lfloor 2^{n/2}t\right\rfloor }\big(f(X_{\left(
j-1\right) 2^{-n/2}})+f(X_{j2^{-n/2}})\big)\left( X_{j2^{-n/2}}-X_{\left( j-1\right) 2^{-n/2}}\right)\\
&=&D_t^{(n)}+E_t^{(n)}.
\end{eqnarray*}%
Since $x^{\kappa }-\mu _{\kappa +1}x$ has centered Hermite rank
$\geq 2$, one can deal with $D_{t}^{\left( n\right) }$ directly
via Theorem \ref{T : BastaJacod}. The conclusion is obtained by
observing that $E_{t}^{\left( n\right) }$ converges u.c.p. to $\mu
_{\kappa +1}\int_{0}^{t}f\left( X_{s}\right) d^\circ X_s$.
\end{proof}

\bigskip

A slight modification of Corollary \ref{cor-odd} yields:

\begin{corollary}\label{cor13}
Let $f:\mathbb{R}\rightarrow\mathbb{R}$ belong to ${\rm C}^2$ with $%
f^{\prime }$ and $f^{\prime \prime }$ bounded, $\kappa \geq 3$ be an odd
integer, and set:
\begin{eqnarray}
\widetilde{J}_{t}^{\left( n\right) }\left( f\right) &=&2^{(\kappa-1)\frac{n}{4}}
\!\!\!\!\sum_{j=1}^{\left%
\lfloor \frac12\big(2^{\frac{n}2}t-1\big)\right\rfloor }\!\!\!\!
f(X_{\left( 2j+1\right) 2^{-\frac{n}2}})\left[
\left( X_{(2j+2)2^{-\frac{n}2}}-X_{\left(2j+1\right)
2^{-\frac{n}2}}\right) ^{\kappa } \right.\nonumber\\
&&\hskip6cm\left.
+\left( X_{(2j+1)2^{-\frac{n}2}}-X_{\left( 2j\right)
2^{-\frac{n}2}}\right) ^{\kappa } \right],\quad t\geq 0.\nonumber
\end{eqnarray}%
Then, as $n\rightarrow +\infty $, the process $\left\{
\widetilde{J}_{t}^{\left( n\right) }\left( f\right)
,X_{t}\right\}_{t\geq 0} $ converges in the sense f.d.d. to
(\ref{st}).
\end{corollary}
\begin{proof}. Follows the proof of Corollary \ref{cor11}.\end{proof}
\bigskip

The next result can be proved analogously.

\bigskip

\begin{corollary}
Let $f:\mathbb{R}\rightarrow\mathbb{R}$ belong to ${\rm C}^2$ with $f'$ and $f''$ bounded, $\kappa \geq
3$ be an odd integer, and define $J_{t}^{\left( n\right) }\left( f\right) $
according to (\ref{Jf-1}). Then, as $n\rightarrow +\infty $, the process $%
\left\{ J_{t}^{\left( n\right) }\left( f\right) ,X_{t}
\right\}_{t\geq 0} $ converges in the sense of f.d.d. to
\begin{equation*}
\left\{
\sqrt{\mu_{2\kappa}}
\int_{0}^{t}f\left( X_{s}\right)dB(s) ,X_{t}
\right\}_{t\geq 0} \text{.}
\end{equation*}
\end{corollary}

\begin{proof}

One can write
\begin{eqnarray*}
J_{t}^{\left( n\right) }\left( f\right) &=&2^{-\frac{n}{4}%
}\sum_{j=1}^{\left\lfloor 2^{n/2}t\right\rfloor }f(X_{\left( j-1\right)
2^{-n/2}})(-1)^j
\left\{ \left( X_{j}^{\left( n\right) }-X_{\left( j-1\right)
}^{\left( n\right) }\right) ^{\kappa }-\mu _{\kappa +1}\left( X_{j}^{\left(
n\right) }-X_{\left( j-1\right) }^{\left( n\right) }\right) \right\} \\
&&+\mu _{\kappa +1}2^{-\frac{n}4}\sum_{j=1}^{\left\lfloor 2^{n/2}t\right\rfloor }f(X_{\left(
j-1\right) 2^{-n/2}})(-1)^j\left( X_{j}^{(n)}-X_{j-1}^{(n)}\right).
\end{eqnarray*}%
Since $x^{\kappa }-\mu _{\kappa +1}x$ has centered Hermite rank $\geq 2$,
Theorem \ref{T :
BastaJacod} gives the desired conclusion.
\end{proof}

\section{Proof of Theorem \ref{thm1}}\label{S : TH1.2}

Fix $t\in[0,1]$, and let, for any $n\in\N$ and $j\in\Z$,
\begin{eqnarray}
U_{j,n}(t)=\sharp\big\{k=0,\ldots,\lfloor 2^nt\rfloor -1 :&&
\label{UPC}
\\ Y(T_{k,n})\!\!\!\!&=&\!\!\!\!j2^{-n/2}\mbox{ and }Y(T_{k+1,n})=(j+1)2^{-n/2}
\big\} \notag\\
D_{j,n}(t)=\sharp\big\{k=0,\ldots,\lfloor 2^nt\rfloor -1:&&
\label{DOWNC}
\\ Y(T_{k,n})\!\!\!\!&=&\!\!\!\!(j+1)2^{-n/2}\mbox{ and }Y(T_{k+1,n})=j2^{-n/2}
\big\}\notag
\end{eqnarray}
denote the number of upcrossings and downcrossings of the interval
$[j2^{-n/2},(j+1)2^{-n/2}]$
within the first $\lfloor 2^n t\rfloor$ steps of the random walk $\{Y(T_{k,n}),\,k\in\N\}$,
respectively.
Also, set
$$
\mathcal{L}_{j,n}(t)=2^{-n/2}\big(U_{j,n}(t)+D_{j,n}(t)\big).
$$

The following statement collects several useful estimates proved
in \cite{KL1}.
\begin{prop}\label{properties-kl}
\begin{enumerate}
\item For every $x\in\R$ and $t\in[0,1]$, we have
$$E\big[| L_t^{x}(Y)|\big] \leq 2\,E\big[|L_1^0(Y)|\big]\,\sqrt{t}\,
{\rm exp}\big(-\frac{x^2}{2t}\big).$$
\item For every fixed $t\in[0,1]$, we have $\sum_{j\in\Z} E\big[|{\cal L}_{j,n}(t)|^2\big]=O(2^{n/2})$.
\item There exists a positive constant $\mu$ such that, for every $a,b\in\R$ with $ab\geq 0$
and $t\in[0,1]$,
$$E\big[|L_t^b(Y)-L_t^a(Y)|^2\big]\leq\mu\,|b-a|\,\sqrt{t}\,
{\rm exp}\left(-\frac{a^2}{2t}\right).$$
\item There exists a random variable $K\in{\rm L}^8$ such that, for every $j\in\Z$, every $n\geq 0$
and every $t\in[0,1]$, one has that
$$|{\cal L}_{j,n}(t)-L_t^{j2^{-n/2}}(Y)|\leq Kn2^{-n/4}\sqrt{L_t^{j2^{-n/2}}(Y)}.$$
\end{enumerate}
\end{prop}
\begin{proof}

The first point is proved in \cite[Lemma 3.3]{KL1}. The proof of
the second point is obtained by simply mimicking the arguments
displayed in the proof of \cite[Lemma 3.7]{KL1}. The third point
corresponds to the content of \cite[Lemma 3.4]{KL1}, while the
fourth point is proved in \cite[Lemma 3.6]{KL1}.
\end{proof}
We will also need the following result:
\begin{proposition}\label{P : Incr2}
Fix some integers $n,N\geq 1$, and let $\kappa$ be an even integer. Then,
as $m\rightarrow\infty $, the random element%
$$
\left\{ X_x,\,2^{-\frac{m}{4}}
\sum_{i=\lfloor \left( j-1\right)2^{\frac{m-n}2}\rfloor +1}^{\lfloor j2^{\frac{m-n}2}\rfloor }
\left[ \left( X_{i}^{\left( m\right)
}-X_{i-1}^{\left( m\right) }\right)^\kappa -\mu_\kappa \right]
\mathcal{L}_{i,m}(t)
:j=-N,\ldots,N
\right\}_{x\in\R,\,t\in[0,1]}
$$
converges weakly in the sense of f.d.d. to
$$
\left\{
X_x,\,\sqrt{\mu_{2\kappa}-\mu_\kappa^2}\int_{(j-1)2^{-n/2}}^{j2^{-n/2}} L_t^x(Y)dB_x :j=-N,\ldots,N
\right\}_{x\in\R,\,t\in[0,1]}.
$$
\end{proposition}
\begin{proof}

For every $m\geq k\geq n$, we can write
$$
2^{-\frac{m}{4}}
\sum_{i=\lfloor \left( j-1\right)2^{\frac{m-n}2}\rfloor +1}^{\lfloor j2^{\frac{m-n}2}\rfloor }
\left[ \left( X_{i}^{\left( m\right)
}-X_{i-1}^{\left( m\right) }\right)^\kappa -\mu_\kappa \right]
\mathcal{L}_{i,m}(t)
=A^{(m,k)}_{j,n,t}+B^{(m,k)}_{j,n,t}+C^{(m,k)}_{j,n,t}
$$
with
\begin{eqnarray*}
A^{(m,k)}_{j,n,t}&=&2^{-\frac{m}{4}}
\sum_{i=\lfloor \left( j-1\right)2^{\frac{k-n}2}\rfloor +1}^{\lfloor j2^{\frac{k-n}2}\rfloor }
L_t^{i2^{-k/2}}(Y)
\sum_{\ell=\lfloor (i-1)2^{\frac{m-k+n}2}\rfloor +1}^{\lfloor i2^{\frac{m-k+n}2}\rfloor }
\left[ \left( X_{\ell}^{\left( m\right)
}-X_{\ell-1}^{\left( m\right) }\right)^\kappa -\mu_\kappa \right]\\
B^{(m,k)}_{j,n,t}&=&2^{-\frac{m}{4}}
\sum_{i=\lfloor \left( j-1\right)2^{\frac{k-n}2}\rfloor +1}^{\lfloor j2^{\frac{k-n}2}\rfloor }
\sum_{\ell=\lfloor (i-1)2^{\frac{m-k}2}\rfloor +1}^{\lfloor i2^{\frac{m-k}2}\rfloor }
\big[L_t^{\ell 2^{-m/2}}(Y)-L_t^{i2^{-k/2}}(Y)\big]
\left[ \left( X_{\ell}^{\left( m\right)
}-X_{\ell-1}^{\left( m\right) }\right)^\kappa -\mu_\kappa \right]\\
C^{(m,k)}_{j,n,t}&=&2^{-\frac{m}{4}}
\sum_{i=\lfloor \left( j-1\right)2^{\frac{m-n}2}\rfloor +1}^{\lfloor j2^{\frac{m-n}2}\rfloor }
\left[ \left( X_{i}^{\left( m\right)
}-X_{i-1}^{\left( m\right) }\right)^\kappa -\mu_\kappa \right]
\big[L_t^{i2^{-m/2}}(Y)-\mathcal{L}_{i,m}(t)\big]
\end{eqnarray*}
We shall study $A^{\left( m,k\right) }$, $B^{\left( m,k\right) }$
and $C^{(m,k)}$
separately. By Proposition \ref{P : Incr}, we know that, as $m\rightarrow
\infty $, the random element%
$$
\left\{ X;2^{-\frac{m}{4}}\!\!\!\!\!\! \sum_{\ell=\lfloor
(i-1)2^{\frac{m-k+n}2}\rfloor +1}^{\lfloor
i2^{\frac{m-k+n}2}\rfloor }\!\! \left[ \left( X_{\ell}^{\left(
m\right) }-X_{\ell-1}^{\left( m\right) }\right)^\kappa -\mu_\kappa
\right] \!:\!-N\le j\le N,\,\lfloor \left(
j-1\right)2^{\frac{k-n}2}\rfloor \! + \! 1\leq i\leq \lfloor
j2^{\frac{k-n}2}\rfloor \right\}
$$
converges in law to
$$
\left\{
X;\sqrt{\mu_{2\kappa}-\mu_\kappa^2}\big( B_{i2^{-k/2}}-B_{(i-1)2^{-k/2}}\big)
:\,-N\le j\le N,\,\lfloor \left( j-1\right)2^{\frac{k-n}2}\rfloor +1\le i
\le \lfloor j2^{\frac{k-n}2\rfloor }
\right\}.
$$
Hence, as $m\rightarrow \infty$, using also the independence between $X$ and $Y$, we have:
$$
\left\{
X;
A_{j,n}^{\left( m,k\right) } :\,
j=- N,\ldots,N
\right\}
\,\,\,\overset{{\rm f.d.d.}}{\Longrightarrow}\,\,\,
\left\{
X;A_{j,n}^{(\infty,k)} :\,
j=- N,\ldots,N
\right\},
$$
where
$$
A_{j,n,t}^{(\infty,k)}\triangleq
\sqrt{\mu_{2\kappa}-\mu_\kappa^2}
\sum_{i=\lfloor \left( j-1\right)2^{\frac{k-n}2}\rfloor +1}^{\lfloor j2^{\frac{k-n}2}\rfloor }
L_t^{i2^{-k/2}}(Y)
\big( B_{i2^{-k/2}}-B_{(i-1)2^{-k/2}}\big).
$$
By letting $k\rightarrow \infty $, one obtains that
$A_{j,n,t}^{(\infty,k)}$ converges in probability towards
$$
\sqrt{\mu_{2\kappa}-\mu_\kappa^2}\int_{(j-1)2^{-n/2}}^{j2^{-n/2}}
L_t^x(Y)dB_x .$$ This proves, by letting $m$ and then $k$ go to
infinity, that $\{X;A_{j,n}^{\left( m,k\right) }:j=-N,\ldots,N\}$
converges in the sense of f.d.d. to
$$
\left\{
X;\sqrt{\mu_{2\kappa}-\mu_\kappa^2}\int_{(j-1)2^{-n/2}}^{j2^{-n/2}} L_t^x(Y)dB_x :j=-N,\ldots,N
\right\}.
$$

To conclude the proof of the Proposition we shall show that, by
letting $m$ and then $k$ go to infinity (for fixed $j$, $n$ and
$t>0$), one has that $\left\vert B_{j,n,t}^{\left( m,k\right)
}\right\vert $ and $\left\vert C_{j,n,t}^{\left( m,k\right)
}\right\vert $ converge to zero in $L^{2}$. Let us first consider
$C_{j,n,t}^{(m,k)}$. In what follows, $c_{j,n}$ denotes a constant
that can be different from line to line. When $t\in[0,1]$ is
fixed, we have, by the independence of Brownian increments and the
first and the fourth points of Proposition \ref{properties-kl}:
\begin{eqnarray*}
E\left[ \left\vert C_{j,n,t}^{\left( m,k\right) }\right\vert^2
\right] &=&2^{-\frac{m}{2}}\,\mathbf{Var}(G^\kappa)
\sum_{i=\lfloor (j-1)2^{\frac{m-n}2}\rfloor +1}^{\lfloor j2^{\frac{m-n}2}\rfloor }
E\left[\left|L_t^{i2^{-m/2}}(Y)-{\cal L}_{i,m}(t)\right|^2\right]\\
&\leq&c_{j,n} \,2^{-m}\,m^2\sum_{i=\lfloor (j-1)2^{\frac{m-n}2}\rfloor +1}^{
\lfloor j2^{\frac{m-n}2}\rfloor }
E\left[\left|L_t^{i2^{-m/2}}(Y)\right|\right]\\
&\leq&c_{j,n} \,2^{-m/2}\,m^2.
\end{eqnarray*}

Let us now consider $B^{(m,k)}_{j,n,t}$. We have,
by the independence of Brownian increments and 
the third point of Proposition \ref{properties-kl}:
\begin{eqnarray*}
E\left[ \left\vert B_{j,n,t}^{\left( m,k\right) }\right\vert^2
\right] &=&2^{-\frac{m}{2}}\,\mathbf{Var}(G^\kappa)
\sum_{i=\lfloor (j-1)2^{\frac{k-n}2}\rfloor+1 }^{\lfloor j2^{\frac{k-n}2}\rfloor +1}
\sum_{\ell=\lfloor (i-1)2^{\frac{m-k}2}\rfloor +1}^{\lfloor i2^{\frac{m-k}2}\rfloor }
 E\left[\left|L_t^{\ell 2^{-m/2}}(Y)-L_t^{i2^{-k/2}}(Y)\right|^2\right]\\
&\leq&c_{j,n} \,2^{-m/2}\sum_{i=\lfloor (j-1)2^{\frac{k-n}2}\rfloor +1}
^{\lfloor j2^{\frac{k-n}2}\rfloor }
\sum_{\ell=\lfloor (i-1)2^{\frac{m-k}2}\rfloor +1}^{\lfloor i2^{\frac{m-k}2}\rfloor }
\big(i2^{-k/2}-\ell\,2^{-m/2}\big)\\
&\leq&c_{j,n} \,2^{-k/2}.
\end{eqnarray*}

The desired conclusion follows immediately.
\end{proof}
The next result will be the key in the proof of the convergence
(\ref{gen2}):
\begin{theorem}
\label{T : BastaKL}
For even $\kappa\geq 2$ and $t\in[0,1]$, set
\begin{eqnarray*}
J_{t}^{\left( n\right) }\left( f\right) &=&2^{-\frac{n}{4}}\,\frac12\sum_{j\in\Z}
\left(
f(X_{\left( j-1\right) 2^{-n/2}})
+f(X_{j\, 2^{-n/2}})\right)
 \left[ 2^{\kappa\frac{n}4}\left( X_{j2^{-\frac{n}{2}}}
-X_{(j-1)2^{-\frac{n}{2}}}\right)^\kappa -\mu_\kappa \right]
\mathcal{L}_{j,n}(t),
\end{eqnarray*}%
where the real-valued function $f$ belong to ${\rm C}^2$ with $f'$ and $f''$ bounded. Then, as $%
n\rightarrow +\infty $, the random element $\left\{ X_x,
J_{t}^{\left( n\right) }\left(
f\right)\right\}_{x\in\R,\,t\in[0,1]}$ converges in the sense of
f.d.d. to
\begin{equation}\label{z2}
\left\{ X_x,
\sqrt{\mu_{2\kappa}-\mu_\kappa^2}\int_{\R}f(X_x)L_t^x(Y)dB_x\right\}_{x\in\R,\,t\in [0,1]}  \text{.}
\end{equation}
\end{theorem}
\begin{proof}

By proceeding as in the beginning of the proof of Theorem \ref{T : BastaJacod}, it is not difficult to
see that the convergence result in the statement is equivalent to the convergence of the pair
$\left\{X_x,\widetilde{J}_t^{(n)}(f)\right\}_{x\in\R,t\in[0,1]}$ to the object in (\ref{z2})
where
$$
\widetilde{J}_t^{(n)}(f)=2^{-\frac{n}{4}}\sum_{j\in\Z}
f(X_{\left( j-1\right) 2^{-n/2}})
 \left[ 2^{\kappa\frac{n}4}\left( X_{j2^{-\frac{n}{2}}}
-X_{(j-1)2^{-\frac{n}{2}}}\right)^\kappa -\mu_\kappa \right]
\mathcal{L}_{j,n}(t).
$$

For every $m\geq n$ and $p\geq 1$, one
has that%
\begin{eqnarray*}
\widetilde{J}_{t}^{\left( m\right) }\left( f\right)  &=&
2^{-\frac{m}{4}}\sum_{j\in\Z} \sum_{i=\lfloor \left(
j-1\right)2^{\frac{m-n}2}\rfloor +1}^{\lfloor
j2^{\frac{m-n}2}\rfloor } \!\!\!\!\!\!f(X_{\left(
i-1\right) 2^{-m/2}}) \left[ \left(X_{i}^{\left( m\right)
}-X_{i-1}^{\left( m\right) }\right)^\kappa -\mu_\kappa \right]
\mathcal{L}_{i,m}(t)
\\
&=&A_{t}^{\left( m,n,p\right) }+B_{t}^{\left( m,n,p\right) }+C_{t}^{\left( m,n,p\right) },
\end{eqnarray*}%
where
\begin{eqnarray*}
A_{t}^{\left( m,n,p\right) } &=&2^{-\frac{m}{4}}\sum_{|i|> p2^{m/2}}
f(X_{\left( i-1\right) 2^{-m/2}})
\left[ \left( X_{i}^{\left( m\right)
}-X_{i-1}^{\left( m\right) }\right)^\kappa -\mu_\kappa \right]
\mathcal{L}_{i,m}(t), \\
B_{t}^{\left( m,n,p\right) }
&=&2^{-\frac{m}{4}}\sum_{|j|\leq p2^{n/2}}
f(X_{\left( j-1\right) 2^{-n/2}})
\sum_{i=\lfloor \left(j-1\right)2^{\frac{m-n}2}\rfloor +1}^{\lfloor j2^{\frac{m-n}2}\rfloor }
\left[ \left( X_{i}^{\left( m\right) }-X_{i-1}^{\left( m\right)
}\right)^\kappa -\mu_\kappa \right]
\mathcal{L}_{i,m}(t),\\
C_{t}^{\left( m,n,p\right) } &=&2^{-\frac{m}{4}}\sum_{|j|\leq
p2^{n/2}} \sum_{i=\lfloor \left( j-1\right)
2^{\frac{m-n}2}\rfloor +1}^{\lfloor j2^{\frac{m-n}2}\rfloor } 
\left[ f(X_{\left( i-1\right)
2^{-m/2}})-f(X_{\left( j-1\right) 2^{-n/2}})
\right] \\
&&\hskip6cm\times
\left[ \left( X_{i}^{\left(
m\right) }-X_{i-1}^{\left( m\right) }\right)^\kappa -\mu_\kappa
\right] \mathcal{L}_{i,m}(t).
\end{eqnarray*}%
We shall study $A^{\left( m,n,p\right) }$, $B^{\left( m,n,p\right) }$
and $C^{(m,n,p)}$
separately. By Proposition \ref{P : Incr2}, we know that, as $m\rightarrow
\infty $, the random element%
$$
\left\{ X;2^{-\frac{m}{4}}
\sum_{i=\lfloor \left( j-1\right)2^{\frac{m-n}2}\rfloor +1}^{\lfloor j2^{\frac{m-n}2}\rfloor }
\left[ \left( X_{i}^{\left( m\right)
}-X_{i-1}^{\left( m\right) }\right)^\kappa -\mu_\kappa \right]
\mathcal{L}_{i,m}(t)
:|j|\leq p2^{n/2}
\right\}
$$
converges in law to
$$
\left\{
X;\sqrt{\mu_{2\kappa}-\mu_\kappa^2}\int_{(j-1)2^{-n/2}}^{j2^{-n/2}} L_t^x(Y)dB_x :|j|\leq p2^{n/2}
\right\}.
$$
Hence, as $m\rightarrow \infty$,
$$
\{X;B^{\left( m,n,p\right) }\}
\,\,\,\overset{{\rm f.d.d.}}{\Longrightarrow}\,\,\,
\{X;B^{\left( \infty,n,p\right) }\}
$$
where
$$
B^{\left( \infty,n,p\right) }_t=
\sqrt{\mu_{2\kappa}-\mu_\kappa^2}
\sum_{|j|\leq p2^{n/2}}
f(X_{\left( j-1\right) 2^{-n/2}})
\int_{(j-1)2^{-n/2}}^{j2^{-n/2}}
L_t^x(Y)dB_x.
$$
By letting $n\rightarrow \infty $, one obtains that
$B^{(\infty,n,p)}$ converges in probability towards
$$B_t^{(\infty,\infty,p)}=\sqrt{\mu_{2\kappa}-\mu_\kappa^2}
\int_{-p}^p f(X_x)L_t^x(Y)dB_x.$$ Finally, by letting
$p\rightarrow\infty$, one obtains, as limit,
$\sqrt{\mu_{2\kappa}-\mu_\kappa^2}\int_{\R} f(X_x)L_t^x(Y)dB_x$.
This proves that, by letting $m$ and then $n$ and finally $p$ go
to infinity, $\{X;B^{\left( m,n,p\right) }\}$ converges in
the sense of f.d.d. to
$\{X;\sqrt{\mu_{2\kappa}-\mu_\kappa^2}\int_{\R}
f(X_x)L_t^x(Y)dB_x\}$.

To conclude the proof of the Theorem we shall show that, by
letting $m$ and then $n$ and finally $p$ go to infinity,
$\left\vert A_{t}^{\left( m,n,p\right) }\right\vert $ and
$\left\vert C_{t}^{\left( m,n,p\right) }\right\vert $ converge to
zero in $L^{2}$.
Let us first consider $C_t^{(m,n,p)}$. When $t\in[0,1]$ is fixed,
the independence of the Brownian increments
yields that%
\begin{eqnarray*}
E\left[ \left\vert C_{t}^{\left( m,n,p\right) }\right\vert^2
\right] &=&2^{-\frac{m}{2}}\,\mathbf{Var}(G^\kappa) \sum_{|j|\leq
p2^{n/2}}\sum_{i=\lfloor \left( j-1\right)
2^{\frac{m-n}2}\rfloor +1}^{\lfloor j2^{\frac{m-n}2}\rfloor } E\left| f(X_{\left(
i-1\right) 2^{-m/2}})-f(X_{\left( j-1\right) 2^{-n/2}})\right|^2\\
&&\hskip7cm\times E\left\{\left|\mathcal{L}_{i,m}(t)\right|^2\right\}\\
&\leq & \mathbf{Var}(G^\kappa) |f'|_\infty^2\, 2^{-\frac{m+n}2}
\sum_{j\in\Z}\sum_{i=\lfloor \left( j-1\right) 2^{\frac{m-n}2}\rfloor +1}
^{\lfloor j2^{\frac{m-n}2}\rfloor }
E\left\{\left|\mathcal{L}_{i,m}(t)\right|^2\right\}\\
&=& \mathbf{Var}(G^\kappa) |f'|_\infty^2\, 2^{-\frac{m+n}{2}} \sum_{i\in\Z}
E\left\{\left|\mathcal{L}_{i,m}(t)\right|^2\right\}.
\end{eqnarray*}%
The second point of Proposition \ref{properties-kl} implies that
$\sum_{i\in\Z}
E\left[\left|\mathcal{L}_{i,m}(t)\right|^2\right]\leq {\rm cst}.\,
2^{m/2}$ uniformly in $t\in[0,1]$. This shows that
\begin{equation*}
\sup_{m,p}\,\,E\left[ \left\vert C_{t}^{\left( m,n,p\right)
}\right\vert^2 \right]\leq {\rm cst.}2^{-n/2}.
\end{equation*}%

Let us now consider $A_t^{(m,n,p)}$. We have
\begin{eqnarray*}
E\left[\vert A_{t}^{\left( m,n,p\right) }\vert^2\right] &= &
\mathbf{Var}(G^\kappa)\,2^{-\frac{m}{2}}\sum_{|i|> p2^{m/2}} E\big|
f(X_{(i-1)2^{-m/2}}) \big|^2
E\left[\vert
\mathcal{L}_{i,m}(t)\vert^2\right]\\
&\leq& \mathbf{Var}(G^\kappa)\,\sup_{t\in[0,1]}E\big|f(X_t)\big|^2\,2^{-\frac{m}{2}}\sum_{|i|> p2^{m/2}}
E\left[\vert \mathcal{L}_{i,m}(t)\vert^2\right].
\end{eqnarray*}
The fourth point in the statement of Proposition
\ref{properties-kl} yields
$$
\vert
\mathcal{L}_{i,m}(t)\vert
\leq L_t^{i2^{-m/2}}(Y)+Km2^{-m/4}\sqrt{L_t^{i2^{-m/2}}(Y)}.
$$
By using the first point in the statement of Proposition
\ref{properties-kl}, we deduce that:
\begin{eqnarray*}
E\left[\vert A_{t}^{\left( m,n,p\right) }\vert^2\right] &\leq&
{\rm cst}.\,2^{-m/2}\,\sum_{i>p2^{m/2}}{\rm exp}\big(-\frac{i^22^{-m}}{2t}\big)\\
&\leq&
{\rm cst}.\,\sum_{i>p2^{m/2}}\int_{(i-1)2^{-m/2}}^{i2^{-m/2}} {\rm exp}\big(-\frac{x^2}{2t}\big)dx\\
&=&
{\rm cst}.\,\int_{p}^{+\infty} {\rm exp}\big(-\frac{x^2}{2t}\big)dx\le \frac{{\rm cst.}}{p}.
\end{eqnarray*}
The desired conclusion follows.
\end{proof}

\noindent We are finally in a position to prove Theorem
\ref{thm1}:

\textbf{Proof of Theorem \ref{thm1}}

{\it Proof of (\ref{gen2})}. By using an equality analogous to
\cite[p. 648, line 8]{KL1} (observe that our definition of
$V_n^{(\kappa)}(f,t)$ is slightly different than the one given in
\cite{KL1}), $2^{(\kappa-3)\frac{n}4}V_n^{(\kappa)}(f,t)$
equals
$$
2^{-\frac{n}{4}}\,\frac12\sum_{j\in\Z}
\left(
f(X_{\left( j-1\right) 2^{-n/2}})
+f(X_{j\, 2^{-n/2}})
\right)
 \left[ 2^{\kappa\frac{n}4}\left( X_{j2^{-n/2}}-
X_{(j-1)2^{-n/2}}\right)^\kappa -\mu_\kappa \right]
\mathcal{L}_{j,n}(t).
$$
As a consequence, (\ref{gen2}) derives immediately from Theorem \ref{T : BastaKL}.\\

{\it Proof of (\ref{gen3})}.
Based on Lemma \ref{lm-kl}, it is showed in \cite{KL1}, p. 658, that
$$
V_n^{(\kappa)}(f,t)=\left\{
\begin{array}{lll}
\frac12\sum_{j=0}^{j^\star-1} \left(
f(X^+_{\left( j-1\right) 2^{-n/2}})
+f(X^+_{j\, 2^{-n/2}})
\right)\left( X^+_{j2^{-n/2}}-
X^+_{(j-1)2^{-n/2}}\right)^\kappa&\mbox{if $j^\star>0$}\\
0&\mbox{if $j^\star=0$}\\
\frac12\sum_{j=0}^{|j^\star|-1} \left(
f(X^-_{\left( j-1\right) 2^{-n/2}})
+f(X^-_{j\, 2^{-n/2}})
\right)
\left( X^-_{j2^{-n/2}}-
X^-_{(j-1)2^{-n/2}}\right)^\kappa&\mbox{if $j^\star<0$}
\end{array}
\right.
$$
Here, $X^+$ (resp. $X^-$) represents $X$ restricted to $[0,\infty)$ (resp. $(-\infty,0]$), and
$j^\star$ is defined as follows:
$$
j^\star=j^\star(n,t)=2^{n/2}\,Y(T_{\lfloor 2^n t\rfloor,n}).
$$
For $t\in[0,1]$, let
$$
Y_n(t)=Y(T_{\lfloor 2^n t\rfloor,n}).
$$
Also, for $t\geq 0$, set
$$
J_n^{\pm}(f,t)=2^{(\kappa-1)\frac{n}{4}}\,\frac12\sum_{j=1}^{\left%
\lfloor 2^{n/2}t\right\rfloor }
\left(
f(X^\pm_{\left( j-1\right) 2^{-n/2}})
+f(X^\pm_{j\, 2^{-n/2}})
\right)
\left( X^{\pm}_{j2^{-n/2}}-X^{\pm}_{\left( j-1\right) 2^{-n/2}}\right)
^{\kappa}
$$
and, for $u\in\R$:
$$
J_n(f,u)=\left\{
\begin{array}{ll}
J_n^{+}(f,u),&\mbox{if $u\geq 0$},\\
J_n^{-}(f,-u),&\mbox{if $u\geq 0$}.\\
\end{array}
\right.
$$
Observe that
\begin{equation}\label{rel-key}
2^{(\kappa-1)\frac{n}4}\,\,V_n^{(\kappa)}(f,t)=J_n\big(f,Y_n(t)\big)\quad\mbox{(see also (4.3) in \cite{KL1})}.
\end{equation}
For every $s,t\in\R$ and $n\geq 1$, we shall prove
\begin{equation}\label{basta-bdg}
E\big|J_n(f,t)-J_n(f,s)\big|^2\leq c_{f,\kappa}
\left(2^{-\frac{n}2}\big| \lfloor 2^{\frac{n}2}t\rfloor - \lfloor
2^{\frac{n}2}s\rfloor\big| + 2^{-n}\big| \lfloor
2^{\frac{n}2}t\rfloor - \lfloor
2^{\frac{n}2}s\rfloor\big|^2\right)
\end{equation}
for a constant $c_{f,\kappa}$ depending only of $f$ and $\kappa$.
For simplicity, we only make the proof when $s,t\geq 0$, but the
other cases can be handled in the same way. For $u\geq 0$, we can
decompose
$$
J_n(f,u)=J_n^{(a)}(f,u)+J_n^{(b)}(f,u)
$$
where
\begin{eqnarray*}
J_n^{(a)}(f,u)&=&2^{(\kappa -1)\frac{n}4}\sum_{j=1}^{\lfloor
2^{\frac{n}2}u\rfloor}
f\big(X_{(j-1)2^{-\frac{n}2}}\big)\big(X_{j2^{-\frac{n}2}}-X_{(j-1)2^{-\frac{n}2}}\big)^\kappa\\
J_n^{(b)}(f,u)&=&\frac12\,2^{(\kappa
-1)\frac{n}4}\sum_{j=1}^{\lfloor 2^{\frac{n}2}u\rfloor}
f'\big(X_{\theta_{j,n}}\big)\big(X_{j2^{-\frac{n}2}}-X_{(j-1)2^{-\frac{n}2}}\big)^{\kappa+1}
\end{eqnarray*}
for some $\theta_{j,n}$ lying between $(j-1)2^{-\frac{n}2}$ and
$j2^{-\frac{n}2}$. By independence, and because $\kappa$ is odd,
we can write, for $0\le s\le t$:
\begin{eqnarray*}
E\big| J_n^{(a)}(f,t)-J_n^{(a)}(f,s)\big|^2
&=&\mu_{2\kappa}\,2^{-\frac{n}2} \sum_{j=\lfloor
2^{\frac{n}2}s\rfloor +1}^{\lfloor 2^{\frac{n}2}t\rfloor}
E\big|f\big(X_{(j-1)2^{-\frac{n}2}}\big)\big|^2\\
&\leq &c_{f,\kappa}\,2^{-\frac{n}2}\big| \lfloor
2^{\frac{n}2}t\rfloor - \lfloor 2^{\frac{n}2}s\rfloor\big|.
\end{eqnarray*}
For $J_n^{(b)}(f,\cdot)$, we have by Cauchy-Schwarz inequality:
$$
E\big| J_n^{(b)}(f,t)-J_n^{(b)}(f,s)\big|^2 \leq
c_{f,\kappa}\,\left(2^{-\frac{n}2}\big| \lfloor
2^{\frac{n}2}t\rfloor - \lfloor
2^{\frac{n}2}s\rfloor\big|\right)^2.
$$
The desired conclusion (\ref{basta-bdg}) follows. Since $X$ and
$Y$ are independent, (\ref{basta-bdg}) yields that
$$E\left|J_n\big(f,Y_n(t)\big)-J_n\big(f,Y(t)\big)\right|^2$$ is
bounded by
\begin{equation*}
c_{f,\kappa} E\left[2^{-\frac{n}2}\big| \lfloor
2^{\frac{n}2}Y_n(t)\rfloor - \lfloor
2^{\frac{n}2}Y(t)\rfloor\big|+ 2^{-n}\big|
\lfloor 2^{\frac{n}2}Y_n(t)\rfloor - \lfloor
2^{\frac{n}2}Y(t)\rfloor\big|^2\right].
\end{equation*}
But this quantity tends to zero as $n\rightarrow\infty$, because
$Y_n(t)\,{\stackrel{{\rm L}^2}{\longrightarrow}}\,Y(t)$ (recall
that $T_{\lfloor 2^n t\rfloor,n}\,{\stackrel{{\rm
L}^2}{\longrightarrow}}\,t$, see Lemma 2.2 in \cite{KL1}).
Combining this latter fact with the independence between
$J_n(f,\cdot)$ and $Y$, and the convergence in the sense of f.d.d.
given by Corollary \ref{cor-odd}, one obtains
$$
J_n(f,\cdot)\,{\longrightarrow}\,\int_0^\cdot
f(X_z)\big(\mu_{\kappa+1}d^\circ X_z+\sqrt{\mu_{2\kappa}-\mu_{\kappa+1}^2}dB_z\big),
$$
where the convergence is in the sense of f.d.d., hence
$$
J_n(f,Y_n)\,{\longrightarrow}\,\int_0^Y
f(X_z)\big(\mu_{\kappa+1}d^\circ X_z+\sqrt{\mu_{2\kappa}-\mu_{\kappa+1}^2}dB_z\big),
$$
where the convergence is once again in the sense of f.d.d.. In
view of (\ref{rel-key}), this concludes the proof of (\ref{gen3}).

\section{Proof of Theorem \ref{thm2}}\label{S : TH1.4}
Let
\begin{eqnarray*}
UU_{2j+1,n}(t)&=&\sharp\big\{k=0,\ldots,\lfloor 2^{n-1}t\rfloor -1:\quad
Y(T_{2k,n})=(2j)2^{-n/2},\\
&&\hskip2cm Y(T_{2k+1,n})=(2j+1)2^{-n/2},\,Y(T_{2k+2,n})=(2j+2)2^{-n/2}\big\}\\
UD_{2j+1,n}(t)&=&\sharp\big\{k=0,\ldots,\lfloor 2^{n-1}t\rfloor -1:\quad
Y(T_{2k,n})=(2j)2^{-n/2},\\
&&\hskip2cm Y(T_{2k+1,n})=(2j+1)2^{-n/2},\,Y(T_{2k+2,n})=(2j)2^{-n/2}\big\}\\
DU_{2j+1,n}(t)&=&\sharp\big\{k=0,\ldots,\lfloor 2^{n-1}t\rfloor -1:\quad
Y(T_{2k,n})=(2j+2)2^{-n/2},\\
&&\hskip2cm Y(T_{2k+1,n})=(2j+1)2^{-n/2},\,Y(T_{2k+2,n})=(2j+2)2^{-n/2}\big\}\\
DD_{2j+1,n}(t)&=&\sharp\big\{k=0,\ldots,\lfloor 2^{n-1}t\rfloor -1:\quad
Y(T_{2k,n})=(2j+2)2^{-n/2},\\
&&\hskip2cm Y(T_{2k+1,n})=(2j+1)2^{-n/2},\,Y(T_{2k+2,n})=(2j)2^{-n/2}\big\}.
\end{eqnarray*}
denote the number of {\sl double} upcrossings and/or downcrossings of the interval
$[(2j)2^{-n/2},(2j+2)2^{-n/2}]$
within the first $\lfloor 2^n t\rfloor$ steps of the random walk $\{Y(T_{k,n}),\,k\in\N\}$.
Observe that
\begin{eqnarray}
S_n^{(\kappa)}(f,t)&=&\sum_{j\in\Z}
f\big(X_{(2j+1)2^{-n/2}}\big)
\left[\big(X_{(2j+2)2^{-n/2}}-X_{(2j+1)2^{-n/2}}\big)^{\kappa}\right.\nonumber\\
&&\hskip1cm \left.+(-1)^{\kappa+1}
\big(X_{(2j+1)2^{-n/2}}-X_{(2j)2^{-n/2}}\big)^{\kappa}\right]
(UU_{2j+1,n}(t)-D\!D_{2j+1,n}(t)\big).\nonumber\\
\label{snk}
\end{eqnarray}
The proof of the following lemma is easily obtained by observing that
the double upcrossings and downcrossings of the interval
$[(2j)2^{-n/2},(2j+2)2^{-n/2}]$ alternate:
\begin{lemme}\label{LemmaUUDD}
Let $t>0$. For each $j\in\Z$,
$$
UU_{2j+1,n}(t)-D\!D_{2j+1,n}(t)=\left\{\begin{array}{lll}
{\bf 1}_{\{0\le j<\widetilde{j^\star}\}}&\mbox{if $\widetilde{j^\star}>0$}\\
0&\mbox{if $\widetilde{j^\star}=0$}\\
-{\bf 1}_{\{\widetilde{j^\star}\le j<0\}}&\mbox{if $\widetilde{j^\star}<0$}
\end{array}\right.
$$
where
$$\widetilde{j^\star}=\widetilde{j^\star}(n,t)=\frac12\,2^{n/2}Y(T_{2\lfloor 2^{n-1}t\rfloor,n}).$$
\end{lemme}
Consequently, by combining Lemma \ref{LemmaUUDD} with (\ref{snk}),
we deduce:
$$
S_n^{(\kappa)}(f,t)=\left\{
\begin{array}{lll}
\sum_{j=0}^{\widetilde{j^\star}-1} f(X^+_{\left( 2j+1\right) 2^{-n/2}})
\left[\left( X^+_{(2j+2)2^{-n/2}}-
X^+_{(2j+1)2^{-n/2}}\right)^\kappa\right.\\
\left.\hskip4cm
+(-1)^{\kappa+1}
\big(X^+_{(2j+1)2^{-n/2}}-X^+_{(2j)2^{-n/2}}\big)^{\kappa}\right]
&\mbox{if $\widetilde{j^\star}>0$}\\
0&\mbox{if $\widetilde{j^\star}=0$}\\
\sum_{j=0}^{|\widetilde{j^\star}|-1} f(X^-_{\left( 2j+1\right) 2^{-n/2}})
\left[\left( X^-_{(2j+2)2^{-n/2}}-
X^-_{(2j+1)2^{-n/2}}\right)^\kappa\right.\\
\left.\hskip4cm
+(-1)^{\kappa+1}
\big(X^-_{(2j+1)2^{-n/2}}-X^-_{(2j)2^{-n/2}}\big)^{\kappa}\right]
&\mbox{if $\widetilde{j^\star}<0$}
\end{array}
\right..
$$
Here, as in the proof of (\ref{gen3}),
$X^+$ (resp. $X^-$) represents $X$ restricted to $[0,\infty)$ (resp. $(-\infty,0]$).
For $t\geq 0$, set
\begin{eqnarray*}
\widetilde{J}_n^{\pm}(f,t)&=&2^{(\kappa-1)\frac{n}{4}}\sum_{j=0}^{
\left\lfloor \frac12\,2^{\frac{n}2}t\right\rfloor
}
f(X^\pm_{\left( 2j+1\right) 2^{-\frac{n}2}})
\left[\left( X^\pm_{(2j+2)2^{-\frac{n}2}}-
X^\pm_{(2j+1)2^{-\frac{n}2}}\right)^\kappa\right.\\
&&\hskip7cm
\left.+(-1)^{\kappa+1}
\big(X^\pm_{(2j+1)2^{-\frac{n}2}}-X^\pm_{(2j)2^{-\frac{n}2}}\big)^{\kappa}\right]
\end{eqnarray*}
and, for $u\in\R$:
$$
\widetilde{J}_n(f,u)=\left\{
\begin{array}{ll}
\widetilde{J}_n^{+}(f,u),&\mbox{if $u\geq 0$},\\
\widetilde{J}_n^{-}(f,-u),&\mbox{if $u\geq 0$}.\\
\end{array}
\right.
$$
Also, let
$$
\widetilde{Y}_n(t)=Y(T_{2\lfloor 2^{n-1}t\rfloor,n}).
$$
Observe that
\begin{equation}\label{rel-key2}
2^{(\kappa-1)\frac{n}4}\,\,S_n^{(\kappa)}(f,t)=\widetilde{J}_n\big(f,\widetilde{Y}_n(t)\big).
\end{equation}
Finally, using Corollary \ref{cor11} (for $\kappa$ even) and
Corollary \ref{cor13} (for $\kappa$ odd), and arguing exactly as
in the proof of (\ref{gen3}), we obtain that the statement of
Theorem \ref{thm2} holds. \fin
\\
{\bf Acknowledgement}. We thank an anonymous referee for insightful remarks.

\end{document}